\RequirePackage[l2tabu,orthodox]{nag}		%
\documentclass{article}
\usepackage{geometry}		%

\usepackage{amsmath}		%
\usepackage{amssymb}		%
\usepackage{amsfonts}		%
\usepackage{amsthm}		%

\usepackage{mathtools}		%

\mathtoolsset{%
}

\usepackage[utf8]{inputenc}		%
\usepackage[T1]{fontenc}		%

\usepackage{dsfont}		%

\usepackage{acronym}		%

\usepackage[labelfont={bf,small},labelsep=colon,font=small]{caption}	%

\usepackage[dvipsnames,svgnames]{xcolor}		%
\colorlet{MyRed}{Crimson!75!black}
\colorlet{MyGreen}{DarkGreen!80!black}
\colorlet{MyBlue}{MediumBlue!80!black}

\usepackage{cancel}		%
\usepackage{latexsym}		%

\usepackage{pifont}		%

\usepackage{subcaption}		%
\usepackage{tikz}		%
\usetikzlibrary{calc,patterns}		%
\usepackage{pgfplots}
\usepgfplotslibrary{groupplots,dateplot}
\usetikzlibrary{patterns,shapes.arrows}
\pgfplotsset{compat=newest}

\usepackage{array}		%
\usepackage{booktabs}		%
\usepackage[inline,shortlabels]{enumitem}		%
\setenumerate{itemsep=\smallskipamount,topsep=\medskipamount}		%

\usepackage{xspace}		%

\usepackage[round]{natbib}		%

\usepackage{hyperref}
\hypersetup{
colorlinks=true,
linktocpage=true,
pdfstartview=FitH,
breaklinks=true,
pdfpagemode=UseNone,
pageanchor=true,
pdfpagemode=UseOutlines,
plainpages=false,
bookmarksnumbered,
bookmarksopen=false,
bookmarksopenlevel=1,
hypertexnames=true,
pdfhighlight=/O,
urlcolor=MyBlue,linkcolor=MyBlue,citecolor=MyBlue,	%
pdftitle={},
pdfauthor={},
pdfsubject={},
pdfkeywords={},
pdfcreator={pdfLaTeX},
pdfproducer={LaTeX with hyperref}
}

\usepackage[sort&compress,capitalize,nameinlink]{cleveref}		%
\crefname{algo}{Algorithm}{Algorithms}
\crefname{assumption}{Assumption}{Assumptions}
\crefformat{equation}{(#2#1#3)}
\crefmultiformat{equation}{(#2#1#3)}%
{ and~(#2#1#3)}{, (#2#1#3)}{ and~(#2#1#3)}
\crefrangeformat{equation}{\upshape(#3#1#4)\textendash(#5#2#6)}

\usepackage{autonum}

\usepackage{algorithm}		%
\usepackage{algpseudocode}		%

\usepackage{thmtools}		%
\usepackage{thm-restate}		%

\theoremstyle{plain}

\newtheorem{theorem}{Theorem}[section]

\newtheorem{corollary}[theorem]{Corollary}		%
\newtheorem{lemma}[theorem]{Lemma}		%

\newtheorem*{corollary*}{Corollary}		%

\theoremstyle{definition}
\newtheorem{assumption}{Assumption}		%

\newtheorem*{definition*}{Definition}		%
\newtheorem*{assumption*}{Assumptions}		%
\newtheorem*{example*}{Example}		%

\theoremstyle{remark}
\newtheorem{remark}{Remark}		%
\newtheorem*{remark*}{Remark}		%

\newcounter{proofpart}

\numberwithin{example}{section}		%

\usepackage[normalem]{ulem}		%

\newcommand{\debug}[1]{#1}		%

\newcommand{\newmacro}[2]{\newcommand{#1}{{\debug{#2}}}}		%
\newcommand{\newop}[2]{\DeclareMathOperator{#1}{\debug{#2}}}		%

\DeclarePairedDelimiter{\bracks}{[}{]}		%
\DeclarePairedDelimiter{\parens}{(}{)}		%

\DeclarePairedDelimiter{\abs}{\lvert}{\rvert}		%

\DeclarePairedDelimiterX{\setdef}[2]{\{}{\}}{#1:#2}		%
\DeclarePairedDelimiterXPP{\exclude}[1]{\mathopen{}\setminus}{\{}{\}}{}{#1}

\newcommand{\R}{\mathbb{R}}		%
\newcommand{\extR}{\R \cup \{+\infty\}}

\DeclareMathOperator*{\argmax}{arg\,max}		%
\DeclareMathOperator{\bigoh}{\mathcal{O}}		%
\DeclarePairedDelimiterXPP{\bigof}[1]{\mathcal{O}}{(}{)}{}{#1}		%

\newcommand{\wstarcl}[1]{\overline{#1}^*}
\DeclareMathOperator{\diam}{diam}		%
\DeclareMathOperator{\dist}{dist}		%
\DeclareMathOperator{\dom}{dom}		%
\DeclareMathOperator{\intr}{int}		%
\DeclareMathOperator{\one}{\mathds{1}}		%
\DeclareMathOperator{\relint}{ri}		%
\DeclareMathOperator{\vol}{vol}		%

\newcommand{\eg}{e.g.,\xspace}		%
\newcommand{\ie}{i.e.,\xspace}		%
\newcommand{\st}{s.t.\xspace}		%
\newcommand{\wrt}{w.r.t.\xspace}		%
\newcommand{\als}{a.s.\xspace}         %

\newcommand{\alt}[1]{#1'}		%
\newcommand{\altalt}[1]{#1''}		%

\newmacro{\dd}{\mathrm{d}}		%
\newcommand{\eps}{\varepsilon}		%

\newmacro{\const}{\alpha}		%
\newmacro{\Const}{\rho}		%
\usepackage{xparse}
\NewDocumentCommand{\coef}{O{\lambda}}{\debug{#1}}
\newmacro{\param}{\theta}		%
\newmacro{\params}{\Theta}		%

\newmacro{\beforestart}{0}		%
\newmacro{\start}{1}		%
\newmacro{\afterstart}{2}		%
\newmacro{\running}{\start,\afterstart,\dotsc}		%
\newmacro{\halfrunning}{1,3/2,2\dotsc}
\newmacro{\half}{\frac{1}{2}}		%

\newmacro{\run}{t}		%
\newmacro{\runalt}{s}		%
\newmacro{\runaltalt}{\tau}		%
\newmacro{\nRuns}{T}		%
\newmacro{\runs}{\mathcal{\nRuns}}		%

\newmacro{\state}{X}		%
\newmacro{\statealt}{Y}		%
\newmacro{\statealtalt}{Z}		%

\newcommand{\init}[1][\state]{\debug{#1}_{\start}}		%

\newmacro{\tstart}{0}		%
\newmacro{\timealt}{s}		%
\newmacro{\horizon}{T}		%

\newmacro{\traj}{x}		%
\newmacro{\trajalt}{y}		%
\newmacro{\trajaltalt}{z}		%

\newmacro{\flowmap}{\Theta}		%
\DeclarePairedDelimiterXPP{\flowof}[2]{\flowmap_{#1}}{(}{)}{}{#2}		%

\newop{\Nash}{NE}		%
\newop{\CE}{CE}		%
\newop{\CCE}{CCE}		%
\newop{\NI}{NI}		%

\newop{\brep}{br}		%
\newop{\preg}{\overline{Reg}}		%
\newop{\val}{val}		%

\newmacro{\play}{i}		%
\newmacro{\playalt}{j}		%
\newmacro{\playaltlalt}{k}		%
\newmacro{\nPlayers}{N}		%
\newmacro{\players}{\mathcal{\nPlayers}}		%

\newmacro{\pure}{\alpha}		%
\newmacro{\purealt}{\beta}		%
\newmacro{\purealtalt}{\gamma}		%
\newmacro{\nPures}{A}		%
\newmacro{\pures}{\mathcal{\nPures}}		%

\newmacro{\loss}{\ell}		%
\newmacro{\pay}{u}		%
\newmacro{\payv}{v}		%
\newmacro{\pot}{f}		%

\newmacro{\game}{\mathcal{G}}		%
\newmacro{\gamefull}{\game(\players,\points,\pay)}		%

\newmacro{\fingame}{\Gamma}		%
\newmacro{\fingamefull}{\Gamma(\players,\pures,\pay)}		%

\newmacro{\gmat}{g}		%
\newmacro{\gdist}{\dist_{\gmat}}
\newmacro{\mfld}{M}		%
\newmacro{\form}{\omega}		%

\newmacro{\tvec}{z}		%
\newmacro{\uvec}{u}		%

\newmacro{\ball}{\overline{\mathbb{B}}}		%
\newmacro{\opball}{\mathbb{B}}
\newmacro{\openball}{\mathbb{B}}		%
\newmacro{\sphere}{\mathbb{S}}		%

\newmacro{\graph}{\mathcal{G}}
\newmacro{\vertices}{\mathcal{V}}
\newmacro{\edges}{\mathcal{E}}

\newmacro{\mat}{M}		%
\newmacro{\hmat}{H}		%

\newmacro{\ones}{\mathbf{1}}		%
\newmacro{\eye}{I}		%
\newmacro{\zer}{\mathbf{0}}		%

\DeclarePairedDelimiter{\norm}{\lVert}{\rVert}		%
\DeclarePairedDelimiterXPP{\dnorm}[1]{}{\lVert}{\rVert}{_{\ast}}{#1}		%

\DeclarePairedDelimiterXPP{\onenorm}[1]{}{\lVert}{\rVert}{_{1}}{#1}		%
\DeclarePairedDelimiterXPP{\twonorm}[1]{}{\lVert}{\rVert}{_{2}}{#1}		%
\DeclarePairedDelimiterXPP{\supnorm}[1]{}{\lVert}{\rVert}{_{\infty}}{#1}		%

\DeclarePairedDelimiterX{\braket}[2]{\langle}{\rangle}{#1,#2}		%
\DeclarePairedDelimiterX{\internalInner}[1]{\langle}{\rangle}{#1}
\usepackage{xparse}
\NewDocumentCommand \inner {s m g}{
\IfBooleanTF{#1}{
\internalInner*{#2\IfValueT{#3}{,#3}}
}
{
\internalInner{#2\IfValueT{#3}{,#3}}
}
}

\newmacro{\vecspace}{\mathcal{V}}		%
\newmacro{\subspace}{\mathcal{W}}		%

\newmacro{\coord}{i}		%
\newmacro{\coordalt}{j}		%
\newmacro{\coordaltalt}{k}		%
\newmacro{\nCoords}{d}		%
\newmacro{\dims}{\nCoords}		%
\newmacro{\vdim}{\nCoords}		%

\newmacro{\pvec}{v}		%

\newmacro{\bvec}{e}		%
\newmacro{\bvecs}{\mathcal{E}}		%

\newmacro{\pspace}{\mathcal{X}}		%
\newmacro{\dspace}{\pspace^{\ast}}		%
\newmacro{\dspacealt}{\mathcal{Y}}		%

\newmacro{\dvec}{\dpoint}		%
\newmacro{\dbvec}{\eps}		%

\newmacro{\dpoint}{x^*}		%
\newmacro{\dpointalt}{y}		%
\newmacro{\dpointaltalt}{\altalt\dpoint}		%
\newmacro{\dpoints}{\mathcal{Y}}		%

\newmacro{\dstate}{Y}		%
\newmacro{\dbase}{w}		%

\newmacro{\dualvar}{\lambda}
\newmacro{\ddualvar}{\lambda^*}
\newmacro{\dualvaralt}{\mu}

\newmacro{\pert}{\Phi}

\newcommand{\defeq}{\coloneqq}		%

\newcommand{\from}{\colon}		%

\newop{\Opt}{Opt}		%
\newop{\Sol}{Sol}		%
\newop{\gap}{Gap}		%
\newop{\orcl}{Or}		%

\newop{\primal}{(P)}
\newop{\dual}{(Q)}

\newmacro{\tfun}{g}		%
\newmacro{\obj}{f}		%
\newmacro{\objalt}{g}		%
\newmacro{\sobj}{F}		%

\newmacro{\gvec}{g}		%
\newmacro{\oper}{A}		%
\newmacro{\vecfield}{v}		%

\newcommand{\sol}[1][\point]{#1^{\star}}		%
\newmacro{\vecsol}{\vecfield(\sol)}		%

\newmacro{\signal}{V}		%
\newmacro{\step}{\gamma}		%
\newmacro{\learn}{\eta}		%

\newmacro{\vbound}{G}		%
\DeclareMathOperator{\lips}{Lip}
\newmacro{\lipcst}{L}
\newmacro{\strong}{\mu}		%
\newmacro{\smooth}{\beta}		%

\newop{\cone}{cone}
\newop{\tspace}{T}		%
\newop{\tcone}{TC}		%
\newop{\dcone}{DC}		%
\newop{\ncone}{NC}		%
\newop{\pcone}{PC}		%
\newop{\hull}{\Delta}		%

\newmacro{\cvx}{\mathcal{C}}		%
\newmacro{\subd}{\partial}		%

\newmacro{\minmax}{\mathcal{L}}		%

\newmacro{\minvar}{\point_{1}}		%
\newmacro{\minvaralt}{\alt\minvar}		%
\newmacro{\minvars}{\points_{1}}		%
\newmacro{\minsol}{\sol[\minvar]}		%

\newmacro{\maxvar}{\point_{2}}		%
\newmacro{\maxvaaltr}{\alt\maxvar}		%
\newmacro{\maxvars}{\points_{2}}		%
\newmacro{\maxsol}{\sol[\maxvar]}		%

\newop{\Eucl}{\Pi}		%
\newop{\logit}{\Lambda}		%
\newop{\dkl}{KL}		%

\newmacro{\hreg}{h}		%
\newmacro{\hconj}{\hreg^{\ast}}		%
\newmacro{\breg}{D}		%
\newmacro{\mprox}{P}		%
\newmacro{\mirror}{Q}		%
\newmacro{\fench}{F}		%
\newmacro{\hstr}{K}		%
\newmacro{\depth}{H}		%
\newmacro{\proxdom}{\points_{\hreg}}		%
\newmacro{\zone}{\mathbb{D}}		%
\newmacro{\bregkernel}{\theta} %

\DeclarePairedDelimiterXPP{\proxof}[2]{\mprox_{#1}}{(}{)}{}{#2}		%

\newmacro{\bregexp}{\alpha}
\newmacro{\bregcst}{M}

\newmacro{\kernelcst}{C}
\newmacro{\kernelexp}{q}

\newmacro{\point}{x}		%
\newmacro{\pointalt}{\alt\point}		%
\newmacro{\pointaltalt}{\altalt\point}		%
\newmacro{\points}{\mathcal{K}}		%
\newmacro{\intpoints}{\relint\points}		%

\newmacro{\basealt}{q}		%
\newmacro{\basealtalt}{u}		%

\newmacro{\open}{\mathcal{U}}		%
\newmacro{\openalt}{\mathcal{V}}		%
\newmacro{\openaltalt}{\mathcal{W}}		%
\newmacro{\closed}{\mathcal{C}}		%
\newmacro{\cpt}{\mathcal{K}}		%
\newmacro{\nbd}{\mathcal{U}}		%
\newmacro{\nhd}{\mathcal{U}}		%
\newmacro{\mset}{A}
\newmacro{\msetalt}{B}

\NewDocumentCommand{\ex}{E{_}{\coupling} o}{\mathbb{E}_{#1}\IfValueT{#2}{\bracks*{#2}}}		%
\newop{\prob}{P}%
\newcommand{\dirac}[1]{\debug{\delta}_{#1}}
\newcommand{\probs}[1][\samples]{\debug{\mathcal{P}}(#1)}
\newop{\empirical}{\prob_{\nsamples}}%
\newmacro{\empex}{\frac{1}{\nsamples}\sum_{\ind=1}^\nsamples}
\newop{\probalt}{Q}
\newop{\coupling}{\pi}
\newop{\couplingalt}{\alt\pi}
\NewDocumentCommand{\couplings}{e{_} O{\samples}}{\debug{\mathcal{P}}\IfValueT{#1}{_{#1}}(#2\times#2)}
\newop{\Var}{Var}		%
\newop{\simplex}{\hull}		%

\newop{\rad}{s}

\DeclarePairedDelimiterXPP{\exof}[1]{\ex}{[}{]}{}{%
 #1}

\DeclarePairedDelimiterXPP{\probof}[2]{#1}{(}{)}{}{%
 #2}

\DeclarePairedDelimiterXPP{\oneof}[1]{\one}{\{}{\}}{}{%
 #1}

\newmacro{\sample}{\xi}		%
\newmacro{\samplealt}{\altsample} %
\newmacro{\altsample}{\zeta}
\newmacro{\samples}{\Xi}		%
\newmacro{\dsample}{v}
\newmacro{\dsamplealt}{w}

\newmacro{\nsamples}{n}
\newmacro{\Nsamples}{N}
\newmacro{\nsamplesalt}{m}

\newmacro{\target}{y}
\newmacro{\Point}{X}
\newmacro{\Target}{Y}
\newmacro{\Pointalt}{\alt\Point}
\newmacro{\Targetalt}{\alt\Target}

\newmacro{\targetalt}{\alt\target}
\newmacro{\targetcost}{\kappa}

\newmacro{\cost}{c}

\NewDocumentCommand{\wass}{s O{\cost} m}{\debug{W}_{#2}\IfBooleanT{#1}{\left}(#3\IfBooleanT{#1}{\right})}

\newmacro{\filter}{\mathcal{F}}		%
\newmacro{\probspace}{(\samples,\filter,\prob)}		%

\newmacro{\radius}{\rho}

\newmacro{\event}{E}       %
\newmacro{\eventalt}{H}       %
\newmacro{\mean}{\mu}		%
\newmacro{\sdev}{\sigma}		%
\newmacro{\variance}{\sdev^{2}}		%

\newmacro{\level}{\alpha}
\newmacro{\levelalt}{p}

\NewDocumentCommand{\Lspace}{E{^}{{1}} O{\prob}}{\debug{L^{#1}(#2)}}

\NewDocumentCommand{\Obj}{e{^} E{_}{\radius} D(){\obj}}{
    {F}\IfValueT{#1}{^{#1}}\IfValueT{#2}{_{#2}}(#3)
}

\NewDocumentCommand{\Objalt}{e{^} E{_}{\radius} D(){\obj}}{
    {G}\IfValueT{#1}{^{#1}}\IfValueT{#2}{_{#2}}(#3)
}

\newmacro{\ind}{i}
\newmacro{\indalt}{j}

\newmacro{\reg}{\regparam}
\newmacro{\regalt}{\regparamalt}
\newmacro{\Reg}{R}
\newmacro{\Regalt}{S}
\newmacro{\Regaltalt}{T}
\newmacro{\regparam}{\varepsilon}
\newmacro{\regparamalt}{\delta}
\newmacro{\regparamaltalt}{\tau}
\newmacro{\Regparam}{\Delta}
\newcommand{\phidiv}[2]{D_\phi(#1 | #2)}

\newmacro{\rv}{X}
\newmacro{\discrprob}{p}
\newmacro{\discrprobalt}{q}

\newmacro{\lebmeas}{Leb}
\newcommand{\lebsp}[1][1]{L^{#1}(\base)}
\newcommand{\lebspb}{L^b(\samples^2)}

\newmacro{\measfuncs}{\lebsp[0]}
\newmacro{\dmeasfuncs}{(\measfuncs)^*}

\NewDocumentCommand{\measures}{O{\samples\times\samples} e{^}}{\mathcal{M}\IfValueT{#2}{^{#2}}(#1)}

\newcommand{\contfuncs}[1][\samples\times\samples]{\mathcal{C}(#1)}
\newcommand{\bdedcontfuncs}[1][\samples\times\samples]{\mathcal{C}_b(#1)}

\newcommand{\bivarfunc}{\varphi}%
\newcommand{\bivarfuncalt}{\psi}
\newcommand{\bivarfuncaltalt}{F}

\newcommand{\tildeobj}{\widetilde{\obj}}

\newcommand{\func}{\objalt}
\newcommand{\funcalt}{h}
\newcommand{\constr}{s}
\newcommand{\constralt}{l}
\newcommand{\funcsup}{\objalt^{\bivarfunc}}

\newcommand{\base}[1][\coupling]{{#1}_0}

\newcommand{\functional}{F}
\newcommand{\functionalalt}{G}

\newcommand{\indicator}[1]{\debug{\iota}_{#1}}

\newmacro{\scalar}{u}

\newmacro{\exponent}{p}

\newmacro{\volcst}{V}

\newmacro{\support}{\sigma}

\newmacro{\modular}{H}
\newcommand{\modspace}[1][]{L^{\modular_{{#1}}}}
\newcommand{\dmodspace}[1][]{(L^\modular_{#1})^*} %
\newcommand{\kdmodspace}[1][]{L^{\modular_{#1}^*}} %
\newmacro{\sigmaalgebra}{\mathcal{F}}
\newmacro{\prodsigmaalgebra}{\sigmaalgebra^{\otimes^2}}

\newcommand{\heartmodspace}[1][]{L^{\modular_{#1}}_{\heartsuit}}

\newcommand{\weaktop}[2]{\sigma(#1, #2)}

\DeclareMathOperator*{\internalesssup}{ess\,sup}
\newcommand{\esssup}[1][\base(\cdot | \sample)]{\internalesssup_{#1}}

\renewcommand{\restriction}{|}

\newmacro{\proper}{\tau}		%

\newmacro{\error}{Z}		%
\newmacro{\noise}{U}		%
\newmacro{\bias}{b}		%
\newmacro{\brown}{W}		%

\newmacro{\serror}{\theta}		%
\newmacro{\snoise}{\xi}		%
\newmacro{\sbias}{\psi}		%

\newmacro{\sbound}{M}		%
\newmacro{\bbound}{B}		%
\newmacro{\noisepar}{\sdev}		%
\newmacro{\noisevar}{\variance}		%

\newcommand{\change}[1]{{#1}}

\newmacro{\fn}{f}		%
\newmacro{\fixmap}{F}		%
\newmacro{\silver}{\Phi}		%
\newmacro{\energy}{E}		%
\newmacro{\Lyap}{L}		%
\newmacro{\diff}{\alpha} 

\newmacro{\thres}{\delta}		%
\newmacro{\basin}{\mathcal{B}}		%
\newmacro{\inhd}{\init[\nbd]}

\newmacro{\seed}{\theta}		%
\newmacro{\seeds}{\Theta}		%
\newmacro{\pdist}{P}		%
\newmacro{\history}{\mathcal{H}}		%

\begin{document}

\title{Regularization for\\ Wasserstein Distributionally Robust Optimization}

\date{}

\newcommand{\address}[2]{%
    \footnote{#2}
    \newcounter{#1}
    \setcounter{#1}{\value{footnote}}%
}

\newcommand{\affiliation}[1]{%
    \footnotemark[\value{#1}]%
}

  \author{%
  Waïss Azizian\address{ENS}{DI, ENS, Univ.~PSL, 75005, Paris, France}\address{UGA}{Univ. Grenoble Alpes, 38000, Grenoble, France}
  \and Franck Iutzeler\affiliation{UGA}%
  \and Jérôme Malick\address{LJK}{CNRS \& LJK, 38000, Grenoble, France}
  }

\maketitle

\begin{abstract}
Optimal transport has recently proved to be a useful tool in various machine learning applications needing comparisons of probability measures. Among these, applications of distributionally robust optimization naturally involve Wasserstein distances in their models of uncertainty, capturing data shifts or worst-case scenarios. 
Inspired by the success of the regularization of Wasserstein distances in optimal transport, we study in this paper the regularization of Wasserstein distributionally robust optimization. First, we derive a general strong duality result of regularized Wasserstein distributionally robust problems. Second, we refine this duality result in the case of entropic regularization and provide an approximation result when the regularization parameters vanish. 
\end{abstract}

\newacro{LHS}{left-hand side}
\newacro{RHS}{right-hand side}
\newacro{iid}[i.i.d.]{independent and identically distributed}
\newacro{lsc}[l.s.c.]{lower semi-continuous}
\newacro{usc}[u.s.c.]{upper semi-continuous}
\newacro{rv}[r.v.]{random variable}
\newacro{NE}{Nash equilibrium}

\newacro{ABP}{abstract Bregman proximal}
\newacro{BP}{Bregman proximal}

\newacro{DGF}{distance-generating function}
\newacro{EG}{extra-gradient}
\newacro{MP}{mirror-prox}
\newacro{MD}{mirror descent}
\newacro{OMD}{optimistic mirror descent}
\newacro{OMWU}{optimistic multiplicative weights update}
\newacro{PMP}{past mirror-prox}
\newacro{AMP}{abstract mirror-prox}
\newacro{MPT}{mirror-prox template}

\newacro{VI}{variational inequality}
\newacro{VIP}{variational inequality problem}
\newacro{KKT}{Karush\textendash Kuhn\textendash Tucker}
\newacro{FOS}{first-order stationary}
\newacro{SOO}{second-order optimality}
\newacro{SOS}{second-order sufficient}
\newacro{DGF}{distance-generating function}
\newacro{SFO}{stochastic first-order oracle}

\newacro{DRO}{distributionally robust optimization}
\newacro{WDRO}{Wasserstein distributionally robust optimization}
\newacro{ML}{machine learning}
\newacro{SVM}{support vector machines}
\newacro{ERM}{empirical risk minimization}
\newacro{OT}{optimal transport}
\newacro{ELBO}{evidence lower bound}
\newacro{MCMC}{Monte Carlo Markov Chain}
\newacro{SAEM}{stochastic approximation expectation-maximization}
\newacro{AD}{automatic differentiation}
\newacro{OR}{operational research}
\newacro{PAC}{probably approximately correct}
\newacro{SA}{stochastic approximation}

\newacro{KL}{Kullback-Leibler}

\section{Introduction} %

\Ac{OT} has a long history and exciting recent developments, notably around applications in machine learning and data science; we refer to the monographs \cite{villaniTopicsOptimalTransportation2003a}, \cite{santambrogioOptimalTransportApplied2015}, \cite{peyre2019computational}, and \cite{merigot2021optimal}. %
One of the key technical properties at the core of recent success of \ac{OT} in these applications is the use of regularization, and specifically entropic regularization, opening the way to efficient computational schemes  (see \eg\cite{cuturiSinkhornDistancesLightspeed2013}) to get theoretically-grounded approximations of the Wasserstein distances.

\Ac{DRO} has recently been formulated using \ac{OT} metrics and has proven to be useful in machine learning (see \eg  \citet{kuhn2019wasserstein}). 
But regularization has still to be studied and used in this context.
In this paper, we propose a study of regularization in \Ac{WDRO}, inspired from several recent developments in \ac{OT}, namely \cite{genevayStochasticOptimizationLargescale2016}, \cite{carlier2017convergence}, \citet{genevay2019sample}, and \cite{patyRegularizedOptimalTransport}. %

\subsection{Distributionally robust optimization with Wasserstein neighborhoods} %

\Ac{DRO} is a popular approach in optimization under uncertainty. %
We briefly present here the ideas and the notation that we will use; we refer to the celebrated paper \cite{esfahani2018data} and the survey paper \cite{kuhn2019wasserstein} for more details and for applications in machine learning.

Standard approaches in stochastic optimization consider 
minimizing the expectation of a random loss with respect to  some input distribution or available data points: 
For an objective $\obj_\param : \samples \to \R$ defined on a sample space $\samples$ and depending on parameters $\param \in \params$, %
this consists in considering
\begin{align}
    \min_{\param \in \params} ~\ex_{\sample \sim \prob}[\obj_\param(\sample)]\,.
\end{align}
Here %
$\prob$ is a fixed probability distribution on $\samples$; in practice, %
it is typically an empirical distribution $\prob = %
\empex \dirac{\sample_\ind}$ stemming from data samples  $(\sample_\ind)_{\ind=1}^\nsamples$.

A \Ac{DRO} counterpart of this problem is to minimize the expectation of the loss with respect to  a set of probability distributions close to $\prob$. 
More precisely, we choose a neighborhood $\nhd(\prob)$ of $\prob$ (called the ambiguity set or the distributional uncertainty region) in the set of probabilities measures on $\samples$, denoted by $\probs$  and consider the worst possible expectation of the objective in this neighborhood. 
The resulting problem is thus of the form
\begin{equation}
    \min_{\param \in \params} ~\sobj(\theta) \qquad \text{where}\quad \sobj(\theta)\defeq \sup_{\probalt \in \nhd(\prob)}\ex_{\sample \sim \probalt}[\obj_\param(\sample)]\,.
    \label{eq:intro_dro}
\end{equation}
A natural way to define the ambiguity set $\nhd(\prob)$ is to consider a ball centered at $\prob$ with radius $\radius > 0$ controlling the required level of robustness. When using the Wassertein distance to define the ball, this gives so-called Wasserstein DRO problems (\ac{WDRO}). %

For a cost function $\cost : \samples\times\samples \to \R_+$,  the Wasserstein distance between two probability distributions $\prob, \probalt\in \probs$ is defined as
\begin{equation}\label{eq:wass}
    \wass[\cost]{\prob, \probalt} = \inf \left\{\ex_{(\sample, \altsample) \sim  \coupling} ~ \cost(\sample, \altsample): \coupling \in \couplings, \coupling_1 = \prob, \coupling_2 = \probalt\right\}\,,
\end{equation}
where $\coupling_1$ and $\coupling_2$ %
denote the first and second marginals of the coupling probability, or transport plan, $\coupling$ defined on $ \samples\times\samples$.
A \ac{WDRO} problem thus has the form \eqref{eq:intro_dro}  with the ambiguity set
\begin{equation}\label{eq:nhd}
    \nhd(\prob) = \left\{\probalt \in \probs: \wass{\prob, \probalt} \leq \radius\right\}\,.
\end{equation}

When the the objective $\obj_\param$ exhibits a simple structure, 
this problem reformulates as a tractable convex optimization problem; see \eg\cite{kuhn2019wasserstein}.
This is exploited in several applications, for instance:
logistic regression (see \eg\cite{yuFastDistributionallyRobust2021}), 
support vector machines (see \eg\cite{shafieezadeh-abadehRegularizationMassTransportation2019}), or
$\ell^1$-regression (see \eg\citep{chenDistributionallyRobustLearning2020}). %
Another argument supporting a \ac{WDRO} approach for machine learning applications
is that it provides generalization guarantees, see \eg \cite{esfahani2018data,an2021generalization}.

\subsection{Contributions, related works, and outline} %

In this paper we study regularization in the context of Wasserstein distributionally robust optimization. 
First, we propose a unified framework for double regularization of the \ac{WDRO} objective function (both in the objective and in the constraint). %
We then provide a strong duality result with general convex regularization functions. This result can be seen as the analogue for \ac{WDRO} of the general result of \cite{patyRegularizedOptimalTransport} for \ac{OT}. %
Second, we refine our analysis in the case of the entropic regularization and obtain an explicit expression for the dual problem. Furthermore, we provide approximation guarantees when the regularization parameters are driven to $0$, \change{adapting the reasoning from \cite{carlier2017convergence}}. These results can be seen as analogues for \ac{WDRO} of results of \citet{genevay2019sample} for the entropic-regularized \ac{OT}. 
\change{Note however, that the reasoning and results used \ac{OT} do not directly apply in the context of \ac{WDRO}. Indeed, in \ac{OT}, the entropic regularization depends on (the product of) the two marginal distributions whereas, in \ac{WDRO}, the second marginal is not fixed but rather an optimization variable, which makes the analysis different and more involved. This distinction is further discussed in the beginning of \cref{sec:entro}.}

\change{Up to our knowledge, regularization in the context of \ac{WDRO} has not been investigated on its own, as we do in this paper.\footnote{\change{Let us mention that studying ``regularization for WDRO'' as we propose here should not be confused with studying ``the regularizing effect of WDRO on learning problems'', which is a separate field of study (see \eg \cite{blanchet2019robust, shafieezadeh-abadehRegularizationMassTransportation2019})}} %
Nevertheless, let us mention the two recent papers related to our developments: \citet{blanchet2020semi} and \cite{wang2021sinkhorn}.
In \citet{blanchet2020semi}, an entropic smoothing of a specific WDRO dual function is introduced and used for computational purposes. Such dual smoothing implicitly corresponds to a regularization of the associated primal problem, but this link is not formally made. 
In contrast, the preprint \cite{wang2021sinkhorn} (which appeared while we were preparing this manuscript) shares similar spirit as our work. 
The entropic regularization of WDRO is proposed and analyzed, with a special focus
on computational aspects. This is complementary to our work which provides a theoretical study of general regularizations as well as approximation guarantees for the entropic regularization. We will come back more precisely on the connections between our results and those of \citet{blanchet2020semi} and \cite{wang2021sinkhorn} in \cref{rem:entrodual}.}

The outline of this paper is the following. The introduction ends below with the definition of the framework of this paper. 
In \cref{sec:reg}, we provide a duality result for a general double regularization, together with an illustration in the case when the transport cost is used a regularization function. In \cref{sec:entro}, we focus our analysis to the entropic regularization to get refined expressions and an explicit control of the quality of the approximation of the underlying \ac{WDRO} problem.

\subsection{Set-up, notation, and assumptions} %

The framework of this paper is the following. 
With $\samples$ a subset of $\R^\dims$, $\prob$ a reference probability distribution over $\samples$, and $\obj \from \samples \to \R$ the underlying objective function (we drop the dependence in $\theta$ to simplify the notation), we consider the sup problem in the objective function of \eqref{eq:intro_dro} with the Wasserstein ball of radius $\rho$ as an ambiguity set\;\eqref{eq:nhd}. Our objective thus writes:
\begin{equation}%
    \sup\setdef{\ex_{\sample \sim \probalt}\obj(\sample)}{\probalt \in \probs,\, \wass{\prob, \probalt } \leq \rho}\,.
\end{equation}
We %
reformulate the above problem, in a concise way, using only couplings as
\begin{equation}\label{eq:wdro_obj}%
 \sup_{\coupling \in \couplings_\prob:~\ex_{\coupling} \cost \leq \radius} \ex_{\coupling_2} \obj \tag{WDRO}
\end{equation}
where $\couplings_{\prob}$ is the set of probability distributions on $\samples\times\samples$ having $\prob$ as a first marginal \begin{equation}
    \couplings_{\prob} \defeq \{\coupling \in \couplings: \coupling_1 = \prob\}.
\end{equation} 
\change{When the space $\samples$ is compact,} we have that 
the topological dual of  $\contfuncs$, the set of continuous functions on\;$\samples\times\samples$, is exactly\;$\measures$, the set of finite signed measures over $\samples\times\samples$ by the Riesz representation \cite[Th.~2.14]{rudinRealComplexAnalysis1987}.
We denote by $\inner{\cdot,\cdot}$ the duality pairing between $\contfuncs[\samples\times\samples]$ and $\measures[\samples\times\samples]$:
    \begin{align}
 \begin{cases}
     \measures[\samples\times\samples] \times \contfuncs[\samples\times\samples] &\longrightarrow \R\\
     (\coupling, \varphi) &\longmapsto \inner{\coupling, \varphi} \defeq \int \varphi\; \dd\coupling\,.
 \end{cases}
    \end{align}

When establishing general duality results, we will also make a constant use of the convex conjugate of a function $\functional\from \contfuncs[\samples\times\samples] \to \extR$ defined as
    \begin{align}
        \functional^*\from
        \begin{cases}
            \measures[\samples\times\samples] &\longrightarrow  \extR\\
            \coupling &\longmapsto \sup_{\bivarfunc \in \contfuncs[\samples]} ~\inner{\coupling, \bivarfunc} - \functional(\bivarfunc)\,,
        \end{cases}
    \end{align}
as well as the preconjugate of a function $\functionalalt\from \measures[\samples\times\samples] \to \extR$ defined as
    \begin{align}
        \functionalalt_*\from
        \begin{cases}
            \contfuncs[\samples\times\samples] &\longrightarrow  \extR\\
            \bivarfunc &\longmapsto \sup_{\coupling \in \measures[\samples]} ~\inner{\coupling, \bivarfunc} - \functionalalt(\coupling)\,.
        \end{cases}
    \end{align}
In presence of convexity, these two operations are dual one another; see \eg \citet[Rem.~5.2]{clasonNonsmoothAnalysisOptimization}. 
More precisely, we have that $(\functional^*)_* = \functional$ when $\functional$ is \ac{lsc}, convex, and proper, and 
$(\functionalalt_*)^* = \functionalalt$ when $\functionalalt$ is weakly-$\star$ \ac{lsc}, convex, and proper. \change{Furthermore, the following duality result will be instrumental in our developments; it is a reformulation \cite[Th.~3.2.6]{botDualityVectorOptimization2009a}, adapted to our purposes.
\begin{lemma}[General Fenchel duality]\label{lemma:fenchel_duality_like}
    Consider a compact subset $\pspace \subset \R^\dims$, a function $h \in \contfuncs[\pspace]$, and two functions $\functional, \functionalalt \from \measures[\pspace] \to \extR$ convex weakly-$\star$ \ac{lsc} proper. If there exists $\bivarfunc \in \dom \functional_* \cap (h - \dom \functionalalt_*)$ such that $\functional_*$ is continuous at $\bivarfunc$,
    then
    \begin{align}
        \sup_{\coupling \in \measures[\pspace]} \inner{\coupling,h} - \functional(\coupling) - \functionalalt(\coupling) = \inf_{\bivarfunc,\bivarfuncalt\in\contfuncs[\pspace]:\bivarfunc+\bivarfuncalt=h} \functional_*(\bivarfunc) + \functionalalt_*(\bivarfuncalt)\,.
    \end{align}
\end{lemma}
 \begin{proof}
    First, the \ac{RHS} can be rewritten as
    \begin{align}
        \inf_{\bivarfunc,\bivarfuncalt\in\contfuncs[\pspace]:\bivarfunc+\bivarfuncalt=h} \functional_*(\bivarfunc) + \functionalalt_*(\bivarfuncalt)
        &=
        \inf_{\bivarfunc \in \contfuncs[\pspace]}  \functional_*(\bivarfunc) + \functionalalt_*(h - \bivarfunc)\,.
    \end{align}
     Then, since $\functional$ is convex weakly-$\star$ \ac{lsc}, we get that $\functional = (\functional_*)^*$ and that $\functional_*$ is proper, convex, and \ac{lsc}; see \citep[Rem.~5.2]{clasonNonsmoothAnalysisOptimization}. The same holds for $\functionalalt$.  
     
     We can thus apply the duality result from \cite[Th.~3.2.6]{botDualityVectorOptimization2009a} with $\functional_*$ and $\functionalalt_*(h-\cdot)$ as primal functions that are proper convex functions, and $\contfuncs[\pspace]$ as primal space. Indeed, the regularity assumption of the lemma exactly gives the regularity condition $(RC_1^{\text{id}})$ of this result. Hence,
    \begin{align}
        \inf_{\bivarfunc \in \contfuncs[\pspace]}  \functional_*(\bivarfunc) + \functionalalt_*(h - \bivarfunc)
        &=
        \sup_{\coupling \in \measures[\pspace]} - (\functional_*)^{*}(-\coupling) - (\functionalalt_*(h-\cdot))^* (\coupling)\\
        &=
        \sup_{\coupling \in \measures[\pspace]} - \functional(-\coupling) - G(-\coupling) - \inner{\coupling,h}\,.
    \end{align}
    Carrying out the change of variable $\coupling \gets - \coupling$ then concludes the proof.
 \end{proof}
 }

\vspace*{1ex}
\section{Regularization of the WDRO objective function}\label{sec:reg}

In this section, we study \eqref{eq:wdro_obj} objectives with additional regularization functions both in the constraints $\ex_{\coupling}\cost \leq \radius$ and in the objective $\ex_{\coupling_2} \obj$. For two arbitrary convex functions $\Reg, \Regalt: \measures \to \extR$, the regularized objective we consider is
\begin{equation}\label{eq:primal}
 \sup_{\coupling \in \couplings_\prob:~\ex_{\coupling} \cost + \Regalt(\coupling) \leq \radius} \ex_{\coupling_2} \obj - \Reg(\coupling) \tag{R-WDRO}.
\end{equation}
We give in \Cref{sec:gene} the expression of the dual of this problem, under some compactness assumptions.  %
We specialize %
in \Cref{sec:ex}
the obtained expression \change{for two specific regularizations, namely the cost regularization and the entropic one. Finally, we present in \Cref{sec:weakening} a general framework, beyond any compactness assumption, where the dual expression still holds.} %

\subsection{Duality for regularized WDRO}\label{sec:gene}

\change{In this section, we give the expression of the dual of \eqref{eq:primal}. This result is the analogue for \ac{WDRO} of the main result of \cite{patyRegularizedOptimalTransport} which gives the dual of the \ac{OT} problem regularized by a general convex function. 
We get this dual expression under the following assumptions, providing a nice duality between objects.
\begin{assumption}\label{ass1}\hfill
    \begin{enumerate}[(i)]
        \item $\samples \subset \R^\dims$ is convex and compact;
        \item $\obj \from \samples \to \R$  and $\cost \from \samples\times\samples \to \R_+$ are both continuous;
        \item For all $\sample$ in $\samples$, $\cost(\sample, \sample) = 0$.
    \end{enumerate}
\end{assumption}}

\change{This formulation of the dual problem can also be seen as a generalization of the existing one for the non-regularized case ($\Reg=\Regalt=0$), see \eg\citet{blanchet2019quantifying,gao2016distributionally}.} \change{Note however that the duality results in these papers rely on weaker assumptions; in particular, $\samples$ is not assumed to be compact and $\obj$, $\cost$ are only upper- and lower-semicontinuous respectively. In the following result, these assumptions are needed to handle general regularizations. We discuss how to alleviate them for particular regularization functions in the next section.}

\begin{theorem}[Strong duality for doubly-regularized WDRO]\label{prop:general_regularized_duality}
    Let \cref{ass1} hold and take two convex, proper, and weakly-$\star$ \ac{lsc} functions $\Reg, \Regalt\colon \measures \to \extR$, such that $\Reg + \Regalt$ is also proper. %
    If the primal problem~\eqref{eq:primal} is strictly feasible (\ie if there exists $\coupling \in \couplings_{\prob} \cap \dom \Reg$ such that $\ex_{\coupling} \cost + \Regalt(\coupling) < \radius$), then we have 
    \begin{equation}\label{eq:DRWDRO}
        \val\eqref{eq:primal} ~=~ \inf_{\dualvar \geq 0} \inf_{\bivarfunc \in \contfuncs} \dualvar \radius + \ex_{\sample \sim \prob} [\sup_{\altsample \in \samples} \obj(\altsample) - \dualvar \cost(\sample, \altsample) - \bivarfunc(\sample, \altsample)] + (\Reg + \dualvar\Regalt)_*(\bivarfunc)\,,
    \end{equation}
    and there exists a primal optimal solution $\sol[\coupling] \in \couplings_{\prob}$ and an optimal dual parameter $\sol[\dualvar] \geq 0$.
\end{theorem}

We prove this theorem by carefully combining two standard duality results (the Lagrangian duality theorem in Banach spaces  \citet{peypouquetConvexOptimizationNormed2015} and the Fenchel duality theorem recalled in \cref{lemma:fenchel_duality_like}) with
a powerful theorem for exchanging minimization/integration \citet{rockafellarVariationalAnalysis1998}. The latter %
is rarely considered in similar contexts; see \citet[Thm.~5]{sinha2018certifying} for an exception.

\begin{proof}%
The first step %
consists in applying the Lagrangian duality theorem of \citet[Thm.~3.68]{peypouquetConvexOptimizationNormed2015}; let us %
check its assumptions. First note that Slater's condition holds by assumption, so we only need to check that the primal problem has at least a solution. This is the case by the following arguments:
(i) the problem is feasible by assumption; (ii) $\couplings$ is weakly-$\star$ sequentially compact (see \eg \citet[Cor.~3.30]{brezisFunctionalAnalysisSobolev2010}); (iii) the constraint set $\{\coupling \in \couplings_\prob:\ex_{\coupling}\cost + \Regalt(\coupling)\leq \radius\}$ is weakly-$\star$ closed (since $\Regalt$ is weakly-$\star$ \ac{lsc} and the constraint $\coupling_1 = \prob$ is weakly-$\star$ closed); and (iv) the objective $\coupling \mapsto \ex_{(\sample,\altsample)\sim\coupling}[\obj(\altsample) - \dualvar \cost(\sample, \altsample)] - \Reg(\coupling)$ is weakly-$\star$ \ac{usc} by assumption. 
As a result, we have Lagrangian duality and existence of a dual solution 
\begin{align}
  \val\eqref{eq:primal}  = \inf_{\dualvar \geq 0}\sup_{\coupling \in \couplings: \coupling_1 = \prob} \ex_{(\sample,\altsample)\sim\coupling} [\obj(\altsample) - \dualvar \cost(\sample, \altsample)] - (\Reg + \dualvar \Regalt)(\coupling) +\dualvar \radius \label{eq:LagDual}
\end{align}

The next step is to write the inner sup %
as an inf. For concision, let us introduce $\Regaltalt_\dualvar \defeq \Reg + \dualvar \Regalt$ and $\bivarfuncaltalt_\dualvar (\sample,\altsample) \defeq \obj(\altsample) - \dualvar \cost(\sample, \altsample)$.
The sup over $\coupling$, for a fixed $\dualvar \geq 0$, thus writes
    \begin{align}
        \sup_{\coupling \in \couplings: \coupling_1 = \prob} \ex_{(\sample,\altsample)\sim\coupling} [\obj(\altsample) - \dualvar \cost(\sample, \altsample)] - \Regaltalt_\dualvar(\coupling) =
        \sup_{\coupling \in\measures} \inner{\coupling, \bivarfuncaltalt_\dualvar} -  \indicator{\couplings_\prob}(\coupling) - \Regaltalt_\dualvar(\coupling)\,\label{eq:LagDual2}
    \end{align}
where $\indicator{\couplings_{\prob}}$ is the indicator function in the sense of convex analysis, \ie for $\coupling \in \measures$,  $\indicator{\couplings_{\prob}}(\coupling) = 0$ if $\coupling \in \couplings$ and $ \coupling_1 = \prob$, and $+\infty$ otherwise. 

\change{Now, we want to apply the duality result of \cref{lemma:fenchel_duality_like} with $\functional \gets \indicator{\couplings_{\prob}}$, $\functionalalt \gets  \Regaltalt_\dualvar$, and $ h\gets\bivarfuncaltalt_\dualvar$. To do so, we need to derive the (pre)conjugates of $\indicator{\couplings_{\prob}}$ and $\Regaltalt_\dualvar$.
By the disintegration theorem, %
any coupling $\coupling(\mathrm{d}\sample, \mathrm{d}\altsample)$ can be written as $\prob(\mathrm{d}\sample) \probalt(\mathrm{d}\altsample | \sample)$ with $\probalt$ a conditional probability on $\samples$. Therefore, 
    \begin{align}\label{eq:proof_extended_duality_formula_first_ineq}
      (\indicator{\couplings_{\prob}})_*(\bivarfunc)   &=\sup\left\{\ex_{(\sample, \altsample) \sim \coupling} \bivarfunc(\sample, \altsample):  \coupling \in \couplings,\, \coupling_1 = \prob\right\}\notag \\
        &= \sup\left\{\ex_{\sample \sim \prob} \ex_{\altsample \sim \probalt(\cdot|\sample)} \bivarfunc(\sample, \altsample): \probalt(\cdot|\cdot)\text{ conditional probability on }\samples\right\}\notag \\
        &\leq \ex_{\sample \sim \prob} [\sup_{\altsample \in \samples} \bivarfunc(\sample, \altsample)] \,.
    \end{align}
    We now proceed to the reverse inequality, noting that a measurable map $\altsample \from \samples \to \samples $ induces a conditional probability $\probalt(\cdot|\sample) = \dirac{\altsample(\sample)}$. Hence, 
    \begin{align}\label{eq:proof_extended_duality_formula_snd_ineq}
       (\indicator{\couplings_{\prob}})_*(\bivarfunc)  &=\sup\left\{\ex_{\sample \sim \prob} \ex_{\altsample \sim \probalt(\cdot|\sample)} \bivarfunc(\sample, \altsample): \probalt(\cdot|\cdot)\text{ conditional probability on }\samples\right\} \notag\\
        &\geq \sup\left\{\ex_{\sample \sim \prob} \bivarfunc(\sample, \altsample(\sample)): \altsample \from \samples \to \samples \text{ measurable}\right\}. 
    \end{align}
    The $\sup$ is finite since $\samples$ is %
    compact and $\bivarfunc$ %
    continuous.
    Let us define $\mathsf{N}\colon\samples \times \R^\dims \longrightarrow\extR$  defined as $\mathsf{N}(\sample, \altsample) =  -\bivarfunc(\sample, \altsample) \text{ if } \altsample \in \samples$ and $+\infty$ otherwise. Since $\bivarfunc$ is continuous and $\samples$ is closed, $\mathsf{N}$ is jointly \ac{lsc}, and, as a consequence, it is a normal integrand by \citep[Ex.~14.31]{rockafellarVariationalAnalysis1998}. So we can apply \citep[Thm.~14.60]{rockafellarVariationalAnalysis1998} to get that 
    \begin{equation}\label{eq:proof_extended_duality_formula}
    \inf\left\{\ex_{\sample \sim \prob} -\bivarfunc(\sample, \altsample(\sample)): \altsample \from \samples \to \samples \text{ measurable} \right\} = 
\ex_{\sample \sim \prob}[ \inf_{\altsample \in \samples}  -\bivarfunc(\sample, \altsample)]\,.
    \end{equation}
    Inverting the signs, we have showed both upper and lower-inequalities, which means that 
\begin{align}
    (\indicator{\couplings_\prob})_*(\bivarfunc) = \ex_{\sample \sim \prob} [\sup_{\altsample \in \samples} \bivarfunc(\sample, \altsample)] \qquad \text{for any $\bivarfunc \in \contfuncs$},
\end{align}
and thus $\dom(\indicator{\couplings_\prob})_* = \contfuncs$. Also, $\bivarfunc \in \contfuncs \mapsto \sup_{\samplealt \in \samples} \bivarfunc(\cdot, \samplealt) \in \contfuncs[\samples]$ is 1-Lipschitz \wrt the norm of the uniform convergence so $(\indicator{\couplings_\prob})_*$ is continuous on its domain.
}

   Finally, since $\Regaltalt_\dualvar$ is convex, proper and weakly-$\star$ \ac{lsc}, $\Regaltalt_{\dualvar*}$ is proper 
    and therefore $\dom \Regaltalt_* \neq \emptyset$. %
    Thus, \cref{lemma:fenchel_duality_like} can be used and gives that
    \begin{align}
        \sup_{\coupling \in \couplings: \coupling_1 = \prob}\!\! \ex_{(\sample,\altsample)\sim\coupling}[\obj(\altsample) - \dualvar \cost(\sample, \altsample)] - \Regaltalt(\coupling)
        =  &\inf_{\bivarfunc \in \contfuncs}  \ex_{\sample \sim \prob} \sup_{\altsample \in \samples} \obj(\altsample) - \dualvar \cost(\sample, \altsample) - \bivarfunc(\sample, \altsample) + \Regaltalt_{\dualvar*}(\bivarfunc)\,, %
    \end{align}
    which, combined with \eqref{eq:LagDual}, leads to the claimed result.
\end{proof}

\vspace*{1ex}
\subsection{\change{Examples of regularized WDRO}}\label{sec:ex}

As an illustration of the duality result, let us consider \change{two cases: first, when the regularization is the transport cost itself and, second, when it is a $\phi$-divergence.}

When the transport cost itself is used a regularization, the expression of the dual simplifies as follow. This expression will be used in the analysis of the next section. %

\begin{corollary}[Duality for cost-regularized \ac{WDRO}]\label{lemma:simple_regularized_dual}
    Let \cref{ass1} hold and take $\regparam,\regparamalt > 0$. Then we have 
    \begin{align}
    \sup_{\coupling \in \couplings_\prob\,:\,\ex_{\coupling} \cost + \regparamalt\, \ex_{\coupling} \cost \leq \radius} \ex_{\coupling_2} \obj - \regparam\, \ex_{\coupling} \cost
    ~=~ \inf_{\dualvar \geq 0}  \dualvar \radius + \ex_{\sample \sim \prob} \sup_{\altsample \in \samples} \left\{\obj(\altsample) - (\regparam + \big(1+\regparamalt)\dualvar\big) \cost(\sample, \altsample)\right\}\,.
    \end{align}
\end{corollary}
\begin{proof}
    Define $\Reg:\coupling \in \measures \mapsto \regparam\inner{\coupling, \cost}$, $\Regalt:\coupling \in \measures \mapsto \regparamalt\inner{\coupling, \cost}$ which are convex, proper and weakly-$\star$ continuous by construction. The preconjugate of their sum is $(\Reg+\dualvar\Regalt)_* = \indicator{\{(\regparam+\dualvar\regparamalt) \cost\}}$.
    Moreover, the primal is strictly feasible thanks to the transport plan $\coupling(\dd\sample, \dd\altsample) =  \prob(\dd\sample) \,\dirac{\sample}(\dd\altsample)$.
    Thus, we can apply \cref{prop:general_regularized_duality} to get the expression.
\end{proof}

\change{
When $\phi$-divergences are used as regularizations, the expression of the dual problem of \eqref{eq:DRWDRO} simplifies as follows. We will come back in more details in \Cref{sec:entro} on the KL-divergence, which is a popular $\phi$-divergence. 
\begin{corollary}[Duality for $\phi$-divergence-regularized \ac{WDRO}]\label{lemma:simple_divergence_regularized_dual}
    Let \cref{ass1} hold and take $\regparam,\regparamalt \geq 0$. Consider $\base \in \couplings$, a convex \ac{lsc} function $\phi : \R \to \extR$ such that $\phi(1) = 0$ and $\phi'(\pm \infty) = \pm \infty$ and define the associated divergence for all $\coupling \in \couplings$ as
    \begin{equation}
        \phidiv{\coupling}{\base} \defeq 
        \begin{cases}
            \int_{\samples^2} \phi \parens*{\frac{\dd \coupling }{\dd \base}} \dd \base &\text{if $\coupling$ is absolutely continuous \wrt $\base$}\\
            +\infty &\text{otherwise}\,.
        \end{cases}
    \end{equation}
    Then, with $\Reg \gets \reg \phidiv{\,\cdot\,}{\base}$ and $\Regalt \gets \regalt \phidiv{\,\cdot\,}{\base}$, if the primal $\eqref{eq:primal}$ is strictly feasible, its value is equal to 
    \begin{align}
    \inf_{\dualvar \geq 0} \inf_{\bivarfuncalt \in \contfuncs} \dualvar \radius + \ex_{\sample \sim \prob} [\sup_{\altsample \in \samples} \obj(\altsample) - \dualvar \cost(\sample, \altsample) - \bivarfuncalt(\sample, \altsample)] + (\regparam + \dualvar \regparamalt)\int_{\samples^2} \phi^* \parens*{\frac{\bivarfuncalt(\sample, \samplealt)}{\regparam + \dualvar \regparamalt}} \dd \base (\sample, \samplealt)\,.
    \end{align}
\end{corollary}
\begin{proof}
    $\Reg$ and $\Regalt$ are convex and proper by construction. 
    Furthermore, $D_\phi$ admits the following variational formula, see \eg \citet[Prop.~4.2.8]{agrawal2021optimal} applied with the pair $(\lebsp[\infty], \measures)$,
    \begin{align}
        \forall \coupling \in \measures,\, 
        \phidiv{\coupling}{\base} = \sup \setdef*{\inner{\coupling}{\bivarfuncalt}  - \int_{\samples^2} \phi^* \circ \bivarfuncalt ~ \dd \base}{\bivarfuncalt \in \lebsp[\infty]}\,.
    \end{align}
    By Lusin's theorem, see \eg \citet[Thm.~2.24]{rudinRealComplexAnalysis1987}, the maximization over $\lebsp[\infty]$ can be replaced by maximization over $\contfuncs$, which guarantees that $D_\phi$ is indeed weak-$\star$ \ac{lsc}.
    Then, we can apply \cref{prop:general_regularized_duality}, using that the preconjugate of $D_\phi$ is exactly $\bivarfuncalt \mapsto \int_{\samples^2} \phi^* \circ \bivarfuncalt ~ \dd \base$ by \citet[Prop.~4.2.6]{agrawal2021optimal}.
\end{proof}
}

\vspace*{1ex}
\subsection{\change{Duality under weaker assumptions}}\label{sec:weakening}

\change{
In this section, we present a generalization of \cref{prop:general_regularized_duality} with weaker assumptions. The main change is that the compactness of $\samples$ is no longer required at the expense of a growth condition on~$\obj$.
Specifically, \Cref{ass1} is replaced by the following set of assumptions.
\begin{assumption}\label{ass2}\hfill
    \begin{enumerate}[(i)]
        \item $\samples \subset \R^\dims$ is a closed %
        set; %
        \item $\obj \from \samples \to \R$  is measurable \wrt the Borel $\sigma$-algebra on $\samples$;
        \item $\cost : \samples \times \samples \to \R_+$ is continuous and, for all $\sample$ in $\samples$, $\cost(\sample, \sample) = 0$;
        \item There is some $\base[\sample] \in \samples$ such that $\obj$ satisfies the growth condition\footnote{\change{In the case when $c$ is a power of a distance $c= d^p$, the growth condition is equivalent to $\abs{\obj(\sample)} \leq  \const \parens{1 + d(\sample, \base[\sample])^p}$, which is standard in unregularized \ac{WDRO}, see \eg \citet{gao2016distributionally,blanchet2020optimal}}}:
            \begin{equation}
                \exists\,\const > 0,\, \text{ \st  }\,
                    \forall \sample, \samplealt \in \samples,\, 
                        \quad\abs{\obj(\samplealt)} \leq  \const \parens{1 + \cost(\sample, \base[\sample]) + \cost(\sample, \samplealt)}\,,
            \end{equation}
        and the integrability condition $\ex_{\sample \sim \prob}[\cost(\sample, \base[\sample])] < + \infty$ holds.
    \end{enumerate}
\end{assumption}

For the sake of readability,  we only state here an informal version of our result, where we omit some technical assumptions on the regularizers and the construction of the spaces involved. The formal statement (\cref{theorem:app_weakened_general_regularized_duality}), along with its proof, is provided in \cref{app:weakening}.

\begin{theorem}[Strong duality for general doubly-regularized WDRO, informal]\label{theorem:weakened_general_regularized_duality}
    Let \cref{ass2} hold and take two convex, proper extended-valued functions $\Reg, \Regalt$ which satisfy some regularity conditions.
    If the primal problem~\eqref{eq:primal} is strictly feasible, then we have 
    \begin{equation}\label{eq:GDRWDRO}
        \val\eqref{eq:primal} = \inf_{\dualvar \geq 0} \inf_{\bivarfunc \in \pspace} \dualvar \radius + \ex_{\sample \sim \prob}\! [\esssup[\altsample \in \samples] \obj(\altsample) - \dualvar \cost(\sample, \altsample) - \bivarfunc(\sample, \altsample)] + (\Reg\restriction_{\dspace} \!+ \dualvar\Regalt\restriction_{\dspace})_*(\bivarfunc)\,,
    \end{equation}
    where $\pspace$ is a Banach function space built from $\cost$, $\Reg$ and $\Regalt$, and the essential supremum is taken \wrt the Lebesgue measure on $\samples$.
\end{theorem}

The gist of the proof consists in carefully crafting $\pspace$ from $\cost$, $\Reg$ and $\Regalt$ and then taking advantage of the duality structure of the pair $(\pspace, \dspace)$. Up to the function spaces involved, the dual expression of \eqref{eq:GDRWDRO} corresponds to the one of 
\eqref{eq:DRWDRO}.
Also, the regularity conditions on $\Reg$ and $\Regalt$ encompass the examples of \cref{sec:ex}, \ie the cost regularization and $\phi$-divergences. Note also that the continuity assumption on the cost function $\cost$ can be removed in the case of $\phi$-divergences; see \cref{app:weakening}. %
}

\vspace*{1ex}
\section{Entropic regularization}\label{sec:entro}

In this section, we specialize and refine the study of the previous section in the case entropic regularization, \ie when the \ac{KL} divergence is used a regularizing function. This regularization is defined for two signed measures with finite variations $\mu, \nu$ as 
\begin{align}
    \dkl(\mu |\nu) =
    \begin{cases}
        \int \log\frac{\dd\mu}{\dd\nu} \,\dd\mu &\text{if $\mu$ and $\nu$ are non-negative and }\mu \ll \nu\\
        +\infty &\text{otherwise}\,
    \end{cases}.
\end{align}
This kind of regularization is very popular in computational \ac{OT}, where it enables the derivation useful approximations of the Wasserstein distance, as, for example the so-called Sinkhorn distance:
\begin{equation}\label{eq:KL}
    \inf_{\coupling \in \couplings:\coupling_1 = \prob,\coupling_2 = \probalt} \ex_{\coupling} \cost + \regparam \dkl(\coupling|\prob \otimes \probalt)\, \qquad \text{for $\prob,\probalt \in \probs$ given}.
\end{equation}
When $\samples$ is compact, the dual of this problem %
is given in \citet[Prop.~2.1]{genevayStochasticOptimizationLargescale2016}
\begin{equation}\label{eq:KLdual}
    \sup_{\obj \in \contfuncs[\samples]} \ex_{\prob} \obj -\regparam \,\ex_{\sample \sim \prob}\log\parens*{\ex_{\altsample \sim \probalt} 
    e^\frac{\obj(\sample) - \cost(\sample,\samplealt)}{\regparam}}\,.
\end{equation}
In \cref{sec:entdual}, we establish a similar result for \ac{WDRO}.
But let us point out here a technical difficulty arising from the \ac{WDRO} framework compared to \ac{OT}. The \ac{KL} divergence \eqref{eq:KL} is taken \wrt the measure $\prob \otimes \probalt$, which does not restrict the set of feasible transport plans (if $\coupling \in \couplings$ satisfies $\coupling_1 = \prob$ and $\coupling_2 = \probalt$, $\coupling$ is indeed absolutely continuous \wrt$\prob \otimes \probalt$). In \ac{WDRO} however, %
only one marginal of the transport plans $\coupling$ is fixed and the support of the optimal coupling can be arbitrary, so that the same regularization as in \eqref{eq:KL} %
cannot directly be used.

We thus propose to regularize \eqref{eq:wdro_obj} with \ac{KL} using a base coupling $\base$ with first marginal $(\base)_1 = \prob$ and consider
\begin{equation}\label{eq:primalentropy}\tag{E-WDRO}
    \sup_{\coupling \in \couplings_\prob:~\ex_{\coupling} \cost + \regparamalt \dkl(\coupling|\base)\leq \radius} \ex_{\coupling_2} \obj - \regparam\dkl(\coupling|\base) .
\end{equation}
The choice of $\base$ restricts the set of %
transport plans 
to those that are absolutely continuous \wrt$\base$. We see in \cref{sec:approx} that this restriction has a very limited impact since a natural choice of\;$\base$ still provides a good approximations of the original problem\;\eqref{eq:wdro_obj} when the regularization parameters $\regparam,\regparamalt$ vanish.

\subsection{Duality for the entropy-regularized problem}\label{sec:entdual}

We derive here a duality theorem for \eqref{eq:primalentropy} involving the \ac{KL} regularization with an arbitrary base coupling $\base \in \couplings_\prob$. The %
result naturally involves the same features as in\;\eqref{eq:KLdual}.

\begin{theorem}[Strong duality for entropy-regularized problems]\label{prop:dual_entropic_regularization}
    Let \cref{ass1} hold, take $\regparam, \regparamalt > 0$, and fix an arbitrary $\base[\coupling] \in \couplings_\prob$. If 
    the primal problem \eqref{eq:primalentropy} is strictly feasible (\ie if there exists $\coupling$ such that $\ex_{\coupling} \cost + \regparamalt \dkl(\coupling|\base) < \radius$), then
    \begin{align}
        \val\eqref{eq:primalentropy} ~=~ \inf_{\dualvar \geq 0} \dualvar \radius + (\regparam + \dualvar \regparamalt) \ex_{\sample \sim \prob} \log\left(\ex_{\altsample \sim \base(\cdot|\sample)} e^{\frac{\obj(\altsample) - \dualvar \cost(\sample, \altsample)}{\regparam + \dualvar \regparamalt}}\right)\,,
        \label{eq:snd_entropic_reg_formula}
    \end{align}
    and there exists a primal optimal solution $\sol[\coupling] \in \couplings_{\prob}$ and an optimal dual parameter $\sol[\dualvar] \geq 0$.
\end{theorem}

The interest of using the entropy as a regularization appears when comparing \eqref{eq:snd_entropic_reg_formula} to the general dual \eqref{eq:DRWDRO}. We see that there is no inf on $\varphi$ and moreover the inner sup is replaced by a smoothed approximation (of the  log-sum-exp type). 

\change{\begin{remark}[Two related results and weaker assumptions]\label{rem:entrodual}
The expression \cref{eq:snd_entropic_reg_formula} already appeared, for special cases, in recent articles/preprints. First,  
the dual problem in \cref{eq:snd_entropic_reg_formula} is proposed in \citet{blanchet2020semi} as a smoothing technique (independently from duality and \ac{KL}-regularization) in a specific context of semi-supervised learning over a finite set $\samples$. %
Second, a very similar duality result (in the case of the KL-regularization in constraints only) appears in \cite[Thm.\;1]{wang2021sinkhorn}. 
We note that this duality result holds under weaker assumptions: their proof is based on arguments specific to the \ac{KL} divergence, whereas the one of \cref{prop:dual_entropic_regularization} is based on \cref{lemma:simple_divergence_regularized_dual} for $\phi$-divergences and thus on  \cref{prop:general_regularized_duality} for generic regularizations. 

It is possible to weaken the assumptions of \cref{prop:dual_entropic_regularization} by using the general \cref{theorem:weakened_general_regularized_duality} instead of \cref{prop:general_regularized_duality}, which essentially replaces \Cref{ass1} by \Cref{ass2}. %
We could also weaken the assumptions by using an alternative proof from duality formulas in variational inference; see \cite[Thm.~2.1]{yoonleeGibbsSampler2021} which is itself inspired by \citet[Lem.~1]{boucheronMomentInequalities2005}. 
We do not provide the details of these two refinements here. Indeed, obtaining duality under weaker assumptions is not what matters most here since the approximation result of the next section requires a combination of compactness and continuity.\qed
\end{remark}}

\begin{proof}
    
   \change{We start by applying \cref{lemma:simple_divergence_regularized_dual} with $\phi: x \in \R_+ \mapsto x \log x -x +1$ (whose conjugate is $\phi^*(y)= e^y - 1$), to get
   }
    \begin{align}
        \val\eqref{eq:primalentropy}=
       \inf_{\dualvar \geq 0} \dualvar \radius +  \underbrace{\inf_{\bivarfunc \in \contfuncs} \ex_{\sample \sim \prob} \left[ \sup_{\altsample \in \samples} \obj(\altsample) - \dualvar \cost(\sample, \altsample) - \bivarfunc(\sample, \altsample) \right] + (\regparam+\dualvar\regparamalt)\inner{\base, e^\frac{\bivarfunc}{\regparam+\dualvar\regparamalt} - 1}}_{\scriptsize(a)}\,.
    \end{align}

For a fixed $\dualvar \geq 0$, we can simplify the expression of the term ${\scriptsize (a)}$ above by introducing $\regparamaltalt \defeq \regparam + \dualvar \regparamalt$, defining $\bivarfuncaltalt_\dualvar \from (\sample,\altsample) \mapsto \obj(\altsample) - \dualvar \cost(\sample, \altsample) \in \contfuncs$, and carrying out the change of variable $\bivarfunc \gets \bivarfuncaltalt_\dualvar - \bivarfunc$. We thus obtain
\begin{align}
{\scriptsize (a)} = &\inf_{\bivarfunc \in \contfuncs} \ex_{\sample \sim \prob} \left[ \sup_{\altsample \in \samples}  \obj(\altsample) - \dualvar \cost(\sample, \altsample) - \bivarfunc(\sample, \altsample)\right] + \regparamaltalt \, \ex_{\base} \left[ e^\frac{\bivarfunc}{\regparamaltalt} - 1 \right] \\
    =
&\inf_{\bivarfunc \in \contfuncs} \ex_{\sample \sim \prob} \left[ \sup_{\altsample \in \samples}  \bivarfunc(\sample, \altsample) \right] + \regparamaltalt \,  \ex_{\base}\left[ e^\frac{\bivarfuncaltalt_\dualvar- \bivarfunc}{\regparamaltalt} - 1\right] .\label{eq:KLReformulationBase}
\end{align}

The rest of the proof is devoted to the reformulation of the term ${\scriptsize (a)}$ expressed as~\cref{eq:KLReformulationBase}.
In order to get rid of the supremum in this expression, %
we restrict the minimization to $\contfuncs[\samples]$ instead of $\contfuncs$\change{, \ie we consider functions $\bivarfunc \in \contfuncs$ of the form $\bivarfunc(\sample, \samplealt) = \func(\sample)$ with $\func \in \contfuncs[\sample]$:}
\begin{align}\label{eq:equality_proof_entropic_reg}
    \val\eqref{eq:KLReformulationBase}\leq
    \inf_{\func \in \contfuncs[\samples]} \ex_{\sample \sim \prob} \left[\func(\sample)\right]  + \regparamaltalt\ex_{(\sample, \altsample)\sim\base} \bracks*{e^\frac{\bivarfuncaltalt_\dualvar(\sample,\altsample)- \func(\sample)}{\regparamaltalt} - 1}.
\end{align}
Let us prove that the inequality above is actually an equality. Fix a bivariate function $\bivarfunc \in \contfuncs$ and consider the univariate function $\funcsup \defeq\sup_{\altsample \in \samples}\bivarfunc(\cdot,\altsample)$. Then, since $- \bivarfunc(\sample, \altsample)\geq -\funcsup(\sample)$, we have
\begin{align}
\label{eq:equality_proof_entropic_reg2}
\ex_{\sample \sim \prob} \left[ \sup_{\altsample \in \samples}  \bivarfunc(\sample, \altsample) \right] + \regparamaltalt   \ex_{\base}\left[ e^\frac{\bivarfuncaltalt_\dualvar- \bivarfunc}{\regparamaltalt} - 1\right]
\geq 
\ex_{\sample \sim \prob} \left[\funcsup(\sample)\right]  + \regparamaltalt\ex_{(\sample, \altsample)\sim\base} \bracks*{e^\frac{\bivarfuncaltalt_\dualvar(\sample,\altsample)- \funcsup(\sample)}{\regparamaltalt} - 1}.
\end{align}
However, $\funcsup$ is not continuous in general but only \ac{lsc}, hence we cannot lower bound the \ac{RHS} of \eqref{eq:equality_proof_entropic_reg2} by the \ac{RHS} of \eqref{eq:equality_proof_entropic_reg}. 
To remedy this issue, we approximate $\funcsup$ with continuous functions defined for $k \geq 1$ by %
$\funcsup_k(\sample) \defeq \inf_{\samplealt \in \samples} \funcsup(\samplealt) + k \norm{\sample - \samplealt}\,$.
Since $\funcsup$ is \ac{lsc}, the functions $(\funcsup_k)$ converge pointwise to $\funcsup$ when $k$ goes to $+\infty$, see \eg \citep[Ex.~9.11]{rockafellarVariationalAnalysis1998}. 
Moreover, these functions are uniformly bounded by $\sup_{\sample \in \samples} |\funcsup(\sample)|$, thus Lebesgue's dominated convergence implies that
\begin{align}
\ex_{\sample \sim \prob} \left[ \sup_{\altsample \in \samples}  \bivarfunc(\sample, \altsample) \right] + \regparamaltalt   \ex_{\base}\left[ e^\frac{\bivarfuncaltalt_\dualvar- \bivarfunc}{\regparamaltalt} - 1\right]
&\geq 
\ex_{\sample \sim \prob} \left[\funcsup(\sample)\right]  + \regparamaltalt\ex_{(\sample, \altsample)\sim\base} \bracks*{e^\frac{\bivarfuncaltalt_\dualvar(\sample,\altsample)- \funcsup(\sample)}{\regparamaltalt} - 1}\\
    &=
\lim_{k \to +\infty}\ex_{\sample \sim \prob} \left[\funcsup_k(\sample)\right]  + \regparamaltalt\ex_{(\sample, \altsample)\sim\base} \bracks*{e^\frac{\bivarfuncaltalt_\dualvar(\sample,\altsample)- \funcsup_k(\sample)}{\regparamaltalt} - 1}\\
    &\geq
  \inf_{\func \in \contfuncs[\samples]} \ex_{\sample \sim \prob} \left[\func(\sample)\right]  + \regparamaltalt\ex_{(\sample, \altsample)\sim\base} \bracks*{e^\frac{\bivarfuncaltalt_\dualvar(\sample,\altsample)- \func(\sample)}{\regparamaltalt} - 1}. \label{eq:equality_proof_entropic_reg3}
\end{align}
Combining \cref{eq:equality_proof_entropic_reg} and \cref{eq:equality_proof_entropic_reg3} gives that
\begin{align}
(a)= &\inf_{\func \in \contfuncs[\samples]} \ex_{\sample \sim \prob} \left[\func(\sample)\right]  + \regparamaltalt\ex_{(\sample, \altsample)\sim\base} \bracks*{e^\frac{\bivarfuncaltalt_\dualvar(\sample,\altsample)- \func(\sample)}{\regparamaltalt} - 1} \label{eq:fst_entropic_reg_formula_just_sup}
\end{align}

The final step of the proof consists in solving the above minimum over $g$. Indeed,  the objective
\begin{align}
\func \mapsto \ex_{\sample \sim \prob} \left[\func(\sample)\right]  + \regparamaltalt\ex_{(\sample, \altsample)\sim\base} \bracks*{e^\frac{\bivarfuncaltalt_\dualvar(\sample,\altsample)- \func(\sample)}{\regparamaltalt} - 1}\,
\end{align}
is convex and differentiable on $\contfuncs[\samples]$. As a consequence, the critical points are minimizers. The gradient at a continuous function  $\func \in \contfuncs[\samples]$ lives in $\measures[\samples]$ and is given by 
\begin{align}
\label{eq:grad_b}
    \prob(\dd \sample) - \parens*{\int_\samples e^\frac{\bivarfuncaltalt_\dualvar(\sample,\altsample)- \func(\sample)}{\regparamaltalt} \base(\dd\samplealt|\sample)}\prob(\dd\sample).
\end{align}
Thus the continuous function
\begin{align}
    \func^\star : \sample \mapsto \regparamaltalt \log \parens*{\int_\samples e^\frac{\bivarfuncaltalt_\dualvar(\sample,\altsample) }{\regparamaltalt} \base(\dd\samplealt|\sample)}
\end{align}
is a solution of~\cref{eq:fst_entropic_reg_formula_just_sup}. We get
\begin{align}
(a) ~=~  \ex_{\sample \sim \prob} \left[\regparamaltalt \log \parens*{\int_\samples e^\frac{\bivarfuncaltalt_\dualvar(\sample,\altsample)}{\regparamaltalt} \base(\dd\samplealt|\sample)}\right]  + \underbrace{ \regparamaltalt \ex_{\sample \sim \prob} \bracks*{\int_\samples e^\frac{\bivarfuncaltalt_\dualvar(\sample,\altsample)- \func^\star(\sample)}{\regparamaltalt} \base(\dd\samplealt|\sample) - 1 } }_{=0 ~~ \text{(by nullity of the gradient in \eqref{eq:grad_b} for $\func^\star$)}}
\end{align}
which in turn gives the desired expression~\cref{eq:snd_entropic_reg_formula}.
\end{proof}

\change{
\begin{remark}[Optimal transport plan]
    In the setting of \cref{prop:dual_entropic_regularization}, the primal optimal solution $\sol[\coupling] \in \couplings_\prob$ can be built explicitly from $\sol[\dualvar]$ by taking for any $\sample \in \samples$
    \begin{equation}
        \sol[\coupling](\dd \samplealt | \sample ) \propto 
        e^\frac{\obj(\samplealt) - \sol[\dualvar] \cost(\sample, \samplealt)}{\reg + \dualvar \regalt} \base(\dd \samplealt | \sample)\,.
    \end{equation}
  Indeed, since the dual function
    \begin{equation}
\dualvar \mapsto \dualvar \radius + (\regparam + \dualvar \regparamalt) \ex_{\sample \sim \prob} \log\left(\ex_{\altsample \sim \base(\cdot|\sample)} e^{\frac{\obj(\altsample) - \dualvar \cost(\sample, \altsample)}{\regparam + \dualvar \regparamalt}}\right)
    \end{equation}
    is differentiable at $\sol[\dualvar]$, the optimality condition for \eqref{eq:snd_entropic_reg_formula} implies that
    $%
        \ex_{\sol[\coupling]}[\cost] + \regalt \dkl \parens*{\sol[\coupling] | \base}
        \leq \radius\,
    $%
     (with equality if and only if $\sol[\dualvar]>0$), 
    \ie that $\sol[\coupling]$ is feasible in \eqref{eq:primalentropy}. Then, one can check that
    \begin{equation}
        \ex_{\coupling_2^*}[\obj] - \reg \dkl \parens*{\sol[\coupling] | \base}
        = 
 \sol[\dualvar] \radius + (\regparam + \sol[\dualvar] \regparamalt) \ex_{\sample \sim \prob} \log\left(\ex_{\altsample \sim \base(\cdot|\sample)} e^{\frac{\obj(\altsample) - \sol[\dualvar] \cost(\sample, \altsample)}{\regparam + \sol[\dualvar] \regparamalt}}\right)
    \end{equation}
    so that, by strong duality (\cref{prop:dual_entropic_regularization}), $\sol[\coupling]$ is optimal for \eqref{eq:primalentropy}.
\end{remark}
}

\subsection{Approximation error of entropy-regularized problems}\label{sec:approx}

In this section, we study the behavior of the approximation error of \eqref{eq:primalentropy} as the regularization parameters vanish to 0.
To quantify this approximation, we specify the cost $c$ and the reference measure $\base$. Specifically, we consider %
that the cost $\cost$ is a norm %
to some power $\exponent \geq 1$ 
 \begin{align}\label{eq:c}
\cost(\sample,\altsample) = \norm{\sample - \altsample}^\exponent
\end{align}
and that the reference measure  $\base\in\couplings_\prob$ is taken, for some $\sdev > 0$, as %
 \begin{align}\label{eq:pi0}
       \base(\dd\sample,\dd\samplealt) ~\propto~ \prob(\dd\sample) \one_{\altsample \in \samples}e^{-\frac{\cost(\sample, \altsample)}{2^{\exponent - 1} \sdev}} \dd\samplealt\,.
    \end{align}
For example, when
$\cost(\sample,\altsample) = \norm{\sample - \altsample}_2^\exponent$ with $\exponent \in \{1,2\}$, $\base(\cdot|\sample)$ is a Laplace or a Gaussian distribution (which is easy to sample from for any $\sample \in \samples$).
We also slightly strengthen \cref{ass1} by assuming that $\samples$ is a convex body and that the functions are Lipschitz continuous.

\begin{theorem}[Approximation for entropic regularization]\label{th:approx-simple}
Let the following conditions hold:
    \begin{enumerate}[\normalfont(i)]
        \item the objective $\obj\from \samples \to \R$ and the cost $\cost\from \samples\times\samples \to \R_+$ are Lipschitz continuous;
        \item the cost $\cost$ 
        and the %
        coupling $\base$ %
        are taken as \eqref{eq:c} and \eqref{eq:pi0}
        with $\sdev > 0$ such that $\ex_{\base} \cost < \radius$;
        \item the set $\samples\subset\R^\dims$ is compact, convex, with nonempty interior.
    \end{enumerate}
    Then, as the regularization parameters $\regparam, \regparamalt > 0$ go to zero, we have
    \begin{align}
        0 \leq
    \val\eqref{eq:wdro_obj} - \val\eqref{eq:primalentropy}
    \leq
    \bigoh\parens*{\dims\,(\regparam + \overline{\dualvar}\regparamalt)\log{\frac{1}{\regparam + \overline{\dualvar}\regparamalt}}}
    \end{align}
  where $\overline{\dualvar} \defeq \frac{2 \sup_{\samples} |\obj|}{\radius - \ex_{\base}\cost}$ is an explicit dual bound.
\end{theorem}

This result is the analogue for \eqref{eq:wdro_obj} of quantitative bounds for \cref{eq:KL} established in \citet{genevay2019sample} in the context of \ac{OT}. The core of the proof consists in introducing a block approximation of the optimal transport plan, following \citet{carlier2017convergence}. In our situation, \change{the second marginal of the transport plan $\coupling$ is not fixed, and it is the variable of the optimization problem. Moreover we have an additional variable $\dualvar$ to take into account. So we have to be careful to modify the approximation scheme of \citet{carlier2017convergence,genevay2019sample} and we introduce an auxiliary regularized problem, that can be better handled that the entropy-regularized problem}. Thus the proof of \cref{th:approx-simple} requires several original steps, as described in the next section.

\subsection{Proof of the approximation theorem}

This section is devoted to the proof of \cref{th:approx-simple}. In fact, we state and prove a slightly more detailed result, formalized in the next theorem, and we show afterwards how Theorem~\ref{th:approx-simple} can be derived as a consequence. The following theorem can indeed be seen as a global version of \cref{th:approx-simple}, with explicit constants and slightly more general assumptions.
To simplify the reading, we denote by the optimal solution of the \emph{entropy}-regularized problem.
    \begin{equation}\label{eq:defF}
    \Obj^{\regparam,\regparamalt} = \sup_{\coupling \in \couplings_\prob:~\ex_{\coupling} \cost + \regparamalt \dkl(\coupling|\base)\leq \radius} \ex_{\coupling_2} \obj - \regparam\dkl(\coupling|\base)\,.
    \end{equation}

\begin{theorem}[Extended approximation theorem]\label{thm:approx_entropic_regularization}
Take a radius $\radius > 0$, regularization parameters $\regparam, \regparamalt > 0$, and suppose that the following conditions hold:
    \begin{enumerate}[\normalfont(i)]
        \item The objective $\obj\from \samples \to \R$ and the cost $\cost\from \samples\times\samples \to \R_+$ are Lipschitz continuous, and that %
        their respective Lipschitz constants satisfy
            $\regparam \leq \lips(\obj) \text{ and } \regparamalt \leq \lips(\cost)\,$;
        \item The cost $\cost$ 
        and the coupling $\base$ %
        are taken as \eqref{eq:c} and \eqref{eq:pi0}
        with $\sdev > 0$ such that $\ex_{\base} \cost < \radius$;    
        \item The set $\samples\subset \R^\dims$  is compact, convex, and satisfies (for $\ball(\sample, \Regparam)$ the %
        ball for $\|\cdot\|$)
            \begin{equation}\label{eq:V}
                \volcst \defeq \inf_{\sample \in \samples,0 < \Regparam \leq d} \frac{\vol\parens*{\samples \cap \ball(\sample, \Regparam)}}{\Regparam^\dims} > 0\,.
            \end{equation}
    \end{enumerate}
    \vspace*{-2ex}
    
    Then, we have,
    \begin{align}
     &   \Obj^{0,0}_{\frac{\radius}{1+\delta/\sdev}}(\obj) -
    (\regparam + \overline{\dualvar}\regparamalt)\parens*{
     \dims
     +   \dims \log\parens*{\frac{\lipcst}{(\regparam + \overline{\dualvar}\regparamalt)\dims}}
    + C
    +\frac{1}{\sdev}\parens*{\frac{(\regparam + \overline{\dualvar}\regparamalt)\dims}{\lipcst}}^\exponent
    } - \frac{\regparam \radius}{\sdev+\delta}\\
    &\leq
        \Obj^{\regparam,\regparamalt}(\obj) \leq
        \Obj^{0,0}(\obj)\,.
    \end{align}
  where $\overline{\dualvar} = \frac{2 \sup_{\samples} |\obj|}{\radius - \ex_{\base}\cost}$, $\lipcst = \lips(\obj) + \overline{\dualvar} \lips(\cost)$, and 
  $C= \min\left\{\log\frac{\vol(\samples)}{\volcst}, \log\frac{I_\sdev}{\volcst}\right\}$ with 
  $I_\sdev= \sdev^{\frac{\dims}{\exponent}}\int_{\R^\dims}e^{-\frac{\norm{\samplealt}^\exponent}{2^{\exponent-1}}}\dd\samplealt$.
\end{theorem}
 
The proof of this result requires a few preliminary steps.
    First, we provide in \cref{lemma:approx_simple_regularization} a simple approximation result for the \emph{cost}-regularized problem
        \begin{equation}\label{eq:defG}
    \Objalt^{\regparam,\regparamalt} = \sup_{\coupling \in \couplings_\prob:~\ex_{\coupling} \cost +\regparamalt \ex_{\coupling} \cost\leq \radius} \ex_{\coupling_2} \obj - \regparam \ex_{\coupling} \cost\,.
    \end{equation}
    Next, we bound in \cref{lem:upper_dual} the dual optimal solution of the \emph{entropy}-regularized problem $\Obj^{\regparam,\regparamalt}$. 
    Finally, for a fixed dual variable, we compare in \cref{lemma:fixedla} the values of the Lagrangians of the \emph{entropy}-regularized problem and the \emph{cost}-regularized one, which is the most technical part of the proof. 
After these three lemmas, we prove the approximation result in the extended version (\cref{thm:approx_entropic_regularization}), and show how the initial version (\cref{th:approx-simple}) can be derived from it.  

\begin{lemma}[Approximation for cost-regularization]\label{lemma:approx_simple_regularization}
    Let \cref{ass1} hold and take $\regparam,\regparamalt > 0$. Then, the following bound hold
    \begin{align}
        \Objalt^{0,0}_{\frac{\radius}{1+\delta}}  - \frac{\regparam \radius}{1+\delta} ~\leq~
        \Objalt^{\regparam, \regparamalt} ~\leq~ \Objalt^{0,0}\,.
    \end{align}
\end{lemma}
\begin{proof}
Since the cost function is non-negative, we directly have $\Objalt^{\regparam, \regparamalt} \leq \Objalt^{0,0}$. From \cref{lemma:simple_regularized_dual}, we write the cost-regularized function $\Objalt^{\regparam,\regparamalt}$ of \eqref{eq:defG} as follows
    \begin{align}
        \Objalt^{\regparam,\regparamalt} %
                     &=   \inf_{\dualvar \geq 0} \dualvar \radius + \ex_{\sample \sim \prob} \left[ \sup_{\altsample \in \samples} \obj(\altsample) - (\regparam + (1+\regparamalt)\dualvar) \cost(\sample, \altsample)\right]\\
                            &=   \inf_{\dualvar' \geq \regparam} \frac{\dualvar' - \regparam}{1 + \regparamalt} \radius + \ex_{\sample \sim \prob}\left[ \sup_{\altsample \in \samples} \obj(\altsample) - \dualvar'\cost(\sample, \altsample)\right] ~~~~~ \text{(with $\dualvar' = \regparam + (1+\regparamalt)\dualvar$)}\\
        &\geq   \inf_{\dualvar' \geq 0}  \frac{\dualvar' - \regparam}{1 + \regparamalt} \radius + \ex_{\sample \sim \prob} \left[ \sup_{\altsample \in \samples} \obj(\altsample) - \dualvar'\cost(\sample, \altsample)\right]\\
        &\geq  \frac{ - \regparam\radius}{1 + \regparamalt} +   \inf_{\dualvar' \geq 0}  \dualvar' \frac{  \radius}{1 + \regparamalt} + \ex_{\sample \sim \prob} \left[ \sup_{\altsample \in \samples} \obj(\altsample) - \dualvar'\cost(\sample, \altsample)\right]\\
        &=
        \sup_{\coupling \in \couplings_\prob:\ex_{\coupling} \cost \leq \frac{\radius}{1+\regparamalt}} \ex_{\coupling_2} \obj - \frac{\regparam \radius}{1+\regparamalt} ~=~ \Objalt^{0,0}_{\frac{\radius}{1+\delta}} - \frac{\regparam \radius}{1+\regparamalt}\,,
    \end{align}
    where the last line follows again from \cref{lemma:simple_regularized_dual} with $\regparam=\regparamalt=0$.
\end{proof}

\begin{lemma}[Upper-bound on dual solutions]\label{lem:upper_dual}
    Under the assumptions of \cref{prop:dual_entropic_regularization}, the optimal solution $\sol[\dualvar]$ of the dual problem of $\Obj^{\regparam,\regparamalt}$ is bounded as follows
\begin{align}\label{eq:proof_approx_entropic_bound_dualvar}
    \sol[\dualvar] \leq \overline{\dualvar} = \frac{2 \sup_{\samples} |\obj|}{\radius - \ex_{\base}\cost}\,.
\end{align}
\end{lemma}

\begin{proof}
\cref{prop:dual_entropic_regularization} gives the existence of $\sol[\dualvar]$, which, by definition, minimizes
\begin{align}
    \func\from \dualvar \mapsto \dualvar \radius + (\regparam+\dualvar\regparamalt) \ex_{\sample \sim \prob} \log\left(\ex_{\altsample \sim \base(\cdot|\sample)} e^{\frac{\obj(\altsample) - \dualvar \cost(\sample, \altsample)}{\regparam+\dualvar\regparamalt}}\right)\,.
\end{align}
On the one hand, $\func(\sol[\dualvar])$ is upper bounded as
\begin{align}
    \func(\sol[\dualvar]) \leq \func(0) &= \regparam\, \ex_{\sample \sim \prob} \log\left(\ex_{\altsample \sim \base(\cdot|\sample)} e^{\frac{\obj(\altsample)}{\regparam}}\right)
    \leq \sup_{\samples}|\obj|\,.
\end{align}
On the other hand, thanks to Jensen's inequality, $\func(\sol[\dualvar])$ is lower-bounded as
\begin{align}
    \func(\sol[\dualvar]) &\geq \sol[\dualvar]\radius + \ex_{(\sample,\altsample) \sim \base}[\obj(\altsample) - \sol[\dualvar] \cost(\sample,\altsample)]                     \geq \sol[\dualvar](\radius - \ex_{\base} \cost) - \sup_{\samples}|\obj|\,.
\end{align}
Combining the two inequalities gives \eqref{eq:proof_approx_entropic_bound_dualvar}.
\end{proof}

\begin{lemma}[Approximation bound for the Lagrangians]\label{lemma:fixedla}
Under the assumptions of \cref{thm:approx_entropic_regularization}, consider 
 \begin{align}
 \Obj^{\regparam,\regparamalt}(\dualvar,\obj)
    = & \sup_{\coupling \in \couplings_\prob} \ex_{(\sample,\altsample)\sim\coupling} [\obj(\altsample) - \dualvar \cost(\sample, \altsample)] - (\regparam + \dualvar \regparamalt)\dkl(\coupling|\base)  , \label{eq:LagDualLemma2}\\
 \text{ and } ~~ \Objalt^{\frac{\regparam}{\sdev},\frac{\regparamalt}{\sdev}}(\dualvar,\obj) 
    = & \sup_{\coupling \in \couplings_\prob} \ex_{(\sample,\altsample)\sim\coupling} [\obj(\altsample) - \dualvar \cost(\sample, \altsample)] - \Big(\frac{\regparam + \dualvar \regparamalt}{\sigma}\Big)\cost(\sample, \altsample). \label{eq:LagDualLemma1}
    \end{align}
Then we have, for a fixed $\Regparam \in (0,\dims]$ and with $I_\sdev(\sample) \defeq \int_{\samples}e^{-\frac{\cost(\sample, \altsample)}{2^{\exponent - 1} \sdev}}\dd \altsample$,
 \begin{align}   
    \Objalt^{\frac{\regparam}{\sdev},\frac{\regparamalt}{\sdev}}(\dualvar, \obj)&\leq
    \Obj^{\regparam,\regparamalt}(\dualvar, \obj)
    + (\lips(\obj) +\!\dualvar \lips(\cost)) \Regparam
    + (\regparam +\dualvar\regparamalt)
    \parens*{
    \frac{\Regparam^\exponent}{\sdev}\!-\! \log\parens*{\volcst\!\Regparam^\dims}
    +
    \ex_{\sample \sim \prob}\log I_\sdev(\sample)\!}.%
    \end{align}
\end{lemma}

\begin{proof}
We start by reformulating $\Objalt^{\frac{\regparam}{\sdev},\frac{\regparamalt}{\sdev}}(\dualvar, \obj)$.
By continuity and compactness, the function $\altsample \mapsto \obj(\altsample) - (\dualvar + (\regparam + \dualvar \regparamalt)/{\sigma}) \cost(\sample,\altsample)$ has a maximizer on $\samples$ for every $\sample$. 
By \citep[Thm.~14.37]{rockafellarVariationalAnalysis1998}, we get that there exists a measurable map $\sol[\altsample]\from \samples \to \samples$ such that $\sol[\altsample](\sample) \in \argmax_{\altsample \in \samples} \obj(\altsample) - (\dualvar + (\regparam + \dualvar \regparamalt)/{\sigma})\cost(\sample,\altsample)$ for any $\sample \in \samples$.
Then,
    \begin{align}
        \sol[\coupling](\dd\sample,\dd\altsample) \defeq \prob(\dd \sample)\, \dirac{\sol[\altsample](\sample)}(\dd\altsample)
    \end{align}
    is an optimal solution, and therefore 
    \begin{align}
    \label{eq:lagcostregopt}
         \Objalt^{\frac{\regparam}{\sdev},\frac{\regparamalt}{\sdev}}(\dualvar, \obj) = \ex_{\sol[\coupling]} \bivarfuncaltalt_\dualvar - \frac{\regparam + \dualvar\regparamalt}{\sdev} \ex_{\sol[\coupling]}\cost.   
    \end{align}
Now %
define $\ball^{\Regparam}(\sample) \defeq \ball(\sol[\altsample](\sample), \Regparam)$ and $\coupling^{\Regparam} \in \couplings_{\prob}$ such that  
\begin{align}
\coupling^\Regparam \propto \one_{\samplealt \in \ball^\Regparam(\sample)}\base(\dd\sample,\dd\samplealt)\,.
\end{align}
Note first that we have
\begin{align}
    \ex_{\sol[\coupling]} \bivarfuncaltalt_\dualvar - \ex_{\coupling^\Regparam} \bivarfuncaltalt_\dualvar
    &= \ex_{\sample \sim \prob}\ex_{\altsample \sim \coupling^\Regparam(\cdot|\sample)}[\bivarfuncaltalt_\dualvar(\sample, \sol[\altsample](\sample)) - \bivarfuncaltalt_\dualvar(\sample, \altsample) ]
    ~\leq~  (\lips(\obj) + \dualvar \lips(\cost)) \Regparam\,, \label{eq:proof_approx_entropic_bound_cost}
\end{align}
since $\bivarfuncaltalt_\dualvar$ is $\big(\lips(\obj) + \dualvar \lips(\cost)\big)$-Lipschitz continuous and the support of $\coupling^\Regparam(\cdot|\sample)$ is $\ball^\Regparam(\sample)$. Now, we proceed to bound $\dkl(\coupling^\Regparam|\base)$ by first noticing that
\begin{align}
    \dkl(\coupling^\Regparam|\base) &= 
    \ex_{\sample \sim \prob}\ex_{\altsample \sim \coupling^\Regparam(\cdot|\sample)} \log\parens*{
        \frac
            {\dd \coupling^\Regparam(\altsample|\sample)}
            {\dd               \base(\altsample|\sample)}
    }%
    \\
    &=
    -\ex_{\sample \sim \prob}\log\parens*{
        \int_{\samples \cap \ball^\Regparam(\sample)}e^{-\frac{\cost(\sample, \altsample)}{2^{\exponent - 1} \sdev}}\dd \altsample
    }
    +
    \ex_{\sample \sim \prob}\log\parens*{
        \int_{\samples}e^{-\frac{\cost(\sample, \altsample)}{2^{\exponent - 1} \sdev}}\dd \altsample \label{eq:KLdelta}
    }
\end{align}

We focus on lower-bounding $\int_{\samples \cap \ball^\Regparam(\sample)}e^{-\frac{\cost(\sample, \altsample)}{2^{\exponent - 1} \sdev}}\dd \altsample$. First, note that by the triangular inequality 
(and since $\exponent\geq1$), for any $\sample, \altsample \in \samples$,
\begin{align}
    \frac{\cost(\sample, \altsample)}{2^\exponent} \leq \frac{\cost(\sample, \sol[\altsample](\sample))+\cost(\sol[\altsample](\sample), \altsample)}{2}\,.
\end{align}
If, in addition, $\altsample$ is in $\ball^\Regparam(\sample)$, this bound becomes $
    \frac{\cost(\sample, \altsample)}{2^{\exponent-1}} \leq \cost(\sample, \sol[\altsample](\sample))+\Regparam^\exponent
$.
Hence, we have
\begin{align}
    \int_{\samples \cap \ball^\Regparam(\sample)}e^{-\frac{\cost(\sample, \altsample)}{2^{\exponent - 1} \sdev}}\dd \altsample
&\geq 
e^{-\frac{\cost(\sample, \sol[\altsample](\sample))+\Regparam^\exponent}{\sdev}} \vol(\samples \cap \ball^{\Regparam}(\sample))\geq 
e^{-\frac{\cost(\sample, \sol[\altsample](\sample))+\Regparam^\exponent}{\sdev}} \volcst \Regparam^\dims\,,
\end{align}
where $\volcst$ is defined in \eqref{eq:V}. 
Plugging the above lower and upper bounds into \eqref{eq:KLdelta} %
yields
\begin{align}
    \dkl(\coupling^\Regparam|\base)
    &\leq
     \frac{\ex_{\sample \sim \prob}\cost(\sample, \sol[\altsample](\sample)) + \Regparam^\exponent}{\sdev}
    - \log\parens*{\volcst \Regparam^\dims}
    + 
    \ex_{\sample \sim \prob}\log I_\sdev(\sample) \\
    &=
    \frac{\ex_{\sol[\coupling]}\cost + \Regparam^\exponent}{\sdev}
    - \log\parens*{\volcst \Regparam^\dims}
    + 
    \ex_{\sample \sim \prob}\log I_\sdev(\sample)\,.\label{eq:proof_approx_entropic_bound_kl} 
\end{align}
Finally, putting \cref{eq:lagcostregopt}, \cref{eq:proof_approx_entropic_bound_cost}, and \cref{eq:proof_approx_entropic_bound_kl} together gives,
\begin{align}
    &\Objalt^{\frac{\regparam}{\sdev},\frac{\regparamalt}{\sdev}}(\dualvar, \obj) = \ex_{\sol[\coupling]} \bivarfuncaltalt_\dualvar - \frac{\regparam+\dualvar\regparamalt}{\sdev} \ex_{\sol[\coupling]}\cost\\
    &= \ex_{\coupling^\Regparam} \bivarfuncaltalt_\dualvar - (\regparam  +\dualvar\regparamalt) \dkl(\coupling^\Regparam|\base) + \parens*{
        \ex_{\sol[\coupling]} \bivarfuncaltalt_\dualvar 
         - \ex_{\coupling^\Regparam} \bivarfuncaltalt_\dualvar }
         + \parens*{ ( \regparam +\dualvar\regparamalt) \dkl(\coupling^\Regparam|\base) 
        - \frac{\regparam+\dualvar\regparamalt}{\sdev} \ex_{\sol[\coupling]}\cost}
    \\
    &\leq
    \Obj^{\regparam,\regparamalt}(\dualvar, \obj)
    + (\lips(\obj) + \dualvar \lips(\cost)) \Regparam
    + (\regparam +\dualvar\regparamalt)
    \parens*{
    \frac{\Regparam^\exponent}{\sdev}
    - \log\parens*{\volcst \Regparam^\dims}
    +
    \ex_{\sample \sim \prob}\log I_\sdev(\sample)
    } %
\end{align}
which is the claimed inequality.
\end{proof}

We have now all the ingredients to establish the extended version of the approximation result.

\begin{proof}[Proof of \cref{thm:approx_entropic_regularization}]
First, notice that by \cref{lemma:approx_simple_regularization}, we have that 
\begin{align}
   \Obj^{0,0}_{\frac{\radius}{1+\delta/\sdev}}  - \frac{\regparam \radius}{\sdev+\delta} =  \Objalt^{0,0}_{\frac{\radius}{1+\delta/\sdev}}  - \frac{\regparam \radius}{\sdev+\delta} \leq
        \Objalt^{\frac{\regparam}{\sdev},\frac{\regparamalt}{\sdev}}.
\end{align}
Thus, using the bound at $\dualvar$ fixed, given by \cref{lemma:fixedla} and the upper-bound \eqref{eq:proof_approx_entropic_bound_dualvar}, we get 
\begin{align}
    &\Objalt^{\frac{\regparam}{\sdev},\frac{\regparamalt}{\sdev}} ~\leq~ \inf_{0 \leq \dualvar \leq \overline{\dualvar}}     \Objalt^{\frac{\regparam}{\sdev},\frac{\regparamalt}{\sdev}}(\dualvar, \obj)\\
    &\leq
    \inf_{0 \leq \dualvar \leq \overline{\dualvar}}
    \Obj^{\regparam,\regparamalt}(\dualvar, \obj)
    + (\lips(\obj) + \dualvar \lips(\cost)) \Regparam
    + (\regparam +\dualvar\regparamalt)
    \parens*{
    \frac{\Regparam^\exponent}{\sdev}
    - \log\parens*{\volcst \Regparam^\dims}
    +
    \ex_{\sample \sim \prob}\log I_\sdev(\sample)
    }\\
    &\leq
    \Obj^{\regparam,\regparamalt}(\obj)
    + \left(\lips(\obj) + \overline{\dualvar} \lips(\cost)\right) \Regparam
    + (\regparam +  \overline{\dualvar}\regparamalt)
    \parens*{
    \frac{\Regparam^\exponent}{\sdev}
    - \log\parens*{\volcst \Regparam^\dims}
    +
    \ex_{\sample \sim \prob}\log I_\sdev(\sample)
    }.%
\end{align}
Minimizing the above bound over $\Regparam > 0$, we get that 
the optimal value $\sol[\Regparam]$
is of the form $\sol[\Regparam] = \frac{(\regparam + \overline{\dualvar}\regparamalt) \dims}{\lipcst} + \bigoh((\regparam + \overline{\dualvar}\regparamalt)^{\exponent + 1})$, \change{with $\lipcst = \lips(\obj) + \overline{\dualvar} \lips(\cost)$, as introduced in the statement of the theorem.}
As a consequence, we set $\Regparam = \frac{(\regparam + \overline{\dualvar}\regparamalt) \dims}{\lipcst}$
in the bound above, which becomes
\begin{align}
    \Objalt^{\frac{\regparam}{\sdev},\frac{\regparamalt}{\sdev}}
    &\leq
    \Obj^{\regparam,\regparamalt}(\obj)
    +(\regparam + \overline{\dualvar}\regparamalt)\parens*{
     \dims
     +   \dims \log\parens*{\frac{\lipcst}{(\regparam + \overline{\dualvar}\regparamalt)\dims}}
    +\ex_{\sample \sim \prob}\log I_\sdev(\sample)
    + \log \frac{1}{\volcst}
    +\frac{1}{\sdev}\parens*{\frac{(\regparam + \overline{\dualvar}\regparamalt)\dims}{\lipcst}}^\exponent
    }.
\end{align}
There is only left to bound the term in $\ex_{\sample \sim \prob}\log I_\sdev(\sample)$.
On one hand, we have\;$e^{-\frac{\cost(\sample, \altsample)}{2^{\exponent - 1} \sdev}}\leq 1$\;so\;that
\[
\int_{\samples}e^{-\frac{\cost(\sample, \altsample)}{2^{\exponent - 1} \sdev}}\dd \altsample \leq  
\vol(\samples).
\]
On the other hand, %
we also have (using the change of variable $\altsample' = {\sdev^{-\frac{1}{\exponent}}}(\altsample - \sample)$ to get $I_\sdev$) 
\[
\int_{\samples}e^{-\frac{\cost(\sample, \altsample)}{2^{\exponent - 1} \sdev}}\dd \altsample 
\leq  %
\int_{\R^\dims}e^{-\frac{\norm{\sample-\samplealt}^\exponent}{2^{\exponent-1}\sdev}}\dd\samplealt
=  I_\sdev\, %
\]
This makes the constant $C$ appear in the bound and thus ends the proof. %
\end{proof}

We finish by explaining how the main theorem, \Cref{th:approx-simple}, stems from \cref{thm:approx_entropic_regularization}. 
On the \ac{LHS} of the inequality in \cref{thm:approx_entropic_regularization}, the unregularized objective has radius $\frac{\radius}{1+\delta/\sigma}$, instead of simply $\radius$ in $\Obj^{0,0}_{{\radius}} = \val\eqref{eq:wdro_obj}$. Thus, we compare in the next lemma the optimal values for these two parameters. %

\begin{lemma}[Comparing optimal values]\label{lem:remark_change_radius}
Under the assumptions of \cref{thm:approx_entropic_regularization}, %
\begin{align}
    \Obj^{0,0}_{{\radius}} \leq \Obj^{0,0}_{\frac{\radius}{1 + \delta/\sigma}} + \bigoh(\delta)\,.
\end{align}
\end{lemma}

\begin{proof}
\change{We first apply \citep[Theorem 5.27]{santambrogioOptimalTransportApplied2015} to get a constant-speed geodesic for the $\exponent$-Wasserstein distance connecting $\prob$ and $\probalt$, which is $\wass{\prob, \probalt}^{\frac{1}{p}}$ with our notation. This means that there exists a family of probability distributions $(\probalt_t)_{t \in [0,1]}$ such that $\probalt_0 = \probalt$, $\probalt_1 = \prob$ and, for any $t \in [0,1]$,
\begin{align}
    \label{eq:ineq_geodesic}
    \wass{\prob, \probalt_t}^{\frac{1}{p}} = (1-t) \wass{\prob, \probalt}^{\frac{1}{p}} \qquad\text{and}\qquad \wass{\probalt_t,\probalt}^{\frac{1}{p}} = t \wass{\prob, \probalt}^\frac{1}{p}\,.
\end{align}
We apply these equations with $\probalt$ such that $\wass{\prob, \probalt} \leq \radius$ and $t = 1 - (1+\regparamalt/\sdev)^{-\frac{1}{p}}$ to obtain
\begin{align}
    \label{eq:ineq_geodesic_2}
    \wass{\prob, \probalt_t} \leq \frac{\radius}{1+\delta/\sigma}\qquad\text{and}\qquad \wass{\probalt, \probalt_t} \leq t^{p}\radius = \bigoh(\regparamalt).
\end{align}
Note that the first inequality above yields
\begin{align}\label{eq:ineq_value_geodesic}
    \ex_{\probalt_t}\obj \leq \sup_{\coupling \in \couplings_\prob:\ex_{\coupling} \cost  \leq \frac{\radius}{1 + \delta/\sigma}} \ex_{\coupling_2} \obj  = \Obj^{0,0}_{\frac{\radius}{1 + \delta/\sigma}}, 
\end{align}
We now use the Kantorovich-Rubinstein inequality (\eg \citet[Thm.~1.14]{villaniTopicsOptimalTransportation2003a}) to write %
\begin{align}
    \ex_{\probalt}\obj - \ex_{\probalt_t} \obj ~\leq~ \wass[1]{\probalt, \probalt_t} \lips(\obj) ~\leq~ \wass{\probalt, \probalt_t}^{\frac{1}{\exponent}} \lips(\obj) \,,
\end{align}
where for the second inequality we used that the $\exponent$-Wasserstein distance is always greater than or equal to the 1-Wasserstein distance (\eg \citet[\S 7.1.2]{villaniTopicsOptimalTransportation2003a}).
Together with~\cref{eq:ineq_value_geodesic}, this %
yields:
\begin{align}
    \ex_{\probalt} \obj \leq \ex_{\probalt_t} \obj + \bigoh(\regparamalt) \leq \Obj^{0,0}_{\frac{\radius}{1 + \delta/\sigma}} + \bigoh(\delta)\,.
\end{align}
Taking the supremum over all %
$\probalt$ such that $\wass{\prob,\probalt}\leq\radius$ allows us to conclude.}
\end{proof}

With the help of the previous lemma, the proof of \Cref{th:approx-simple} comes easily from \Cref{thm:approx_entropic_regularization}.

\begin{proof}[Proof of \Cref{{th:approx-simple}}]
We start with checking that the fact that $\samples$ is a compact convex body (condition (iv) in \Cref{th:approx-simple}) implies that $\volcst >0$  (condition (iv) in \Cref{thm:approx_entropic_regularization}). Introduce, for $\sample \in \samples$, the function
\[
\nu_{\xi}: \Regparam \mapsto \frac{\vol\parens*{\samples \cap \ball(\sample, \Regparam)}}{\Regparam^\dims} = {\vol\parens*{\frac{1}{{\Regparam}}(\samples - {\sample})\cap \ball(0,1)}}.
\]
Since $\samples$ is convex, we easily get that $\nu_{\xi}$ is non-increasing.
Thus we can lower-bound the constant $\volcst$ as follows, using $\diam(\samples) \defeq \sup_{\sample, \samplealt \in \samples} \norm{\sample - \samplealt}$ the diameter of $\samples$.
\begin{align}
    \volcst &= %
    \inf_{\sample \in \samples,0 < \Regparam \leq d} \nu_{\xi}(\Regparam) ~\geq ~ \inf_{\sample \in \samples,0 < \Regparam \leq \max(\dims, \diam(\samples))} \nu_{\xi}(\Regparam)\\
            &= \inf_{\sample \in \samples} \frac{\vol\parens*{\samples \cap \ball(\sample, \max(\dims, \diam(\samples))}}{(\max(\dims, \diam(\samples))^\dims} ~\geq~ \frac{\vol\parens*{\samples}}{(\max(\dims, \diam(\samples))^\dims} > 0.
\end{align}
So we get $\volcst >0$ which is condition (iv) in \Cref{thm:approx_entropic_regularization}.
Thus we can apply \Cref{thm:approx_entropic_regularization} %
and \Cref{lem:remark_change_radius}.
Noting that $\Obj^{\regparam,\regparamalt}(\obj) = \val\eqref{eq:primalentropy}$ and $\Obj^{0,0}_{{\radius}} = \val\eqref{eq:wdro_obj}$, and 
combining the obtained bound with~\Cref{lem:remark_change_radius} gives the result.
\end{proof}

\section{Conclusion, perspectives}
Inspired by the success of regularization in \ac{OT}, we proposed and studied %
a regularization scheme for \ac{WDRO} problems. We derived the expression of the dual objective function in the general case as well as a refined one in the particular setting of the entropic regularization. In addition, we showed that the difference between the original \ac{WDRO} problem and the entropic one is properly controlled by the regularization parameters.

\change{Since regularization in \ac{OT} has shown attractive computational advantages and %
statistical benefits, an exciting research direction is to investigate whether similar gains hold for regularization in~\ac{WDRO}. It is also worth studying possible extensions of the results of this paper, including the approximation errors with general regularizations, beyond the entropic case.}

\section*{Acknowledgments}
We thank the associate editor and the two referees for their careful reading and numerous suggestions. In particular, a referee suggested us to study duality beyond compactness/continuity, which has led to the developments sketched in \cref{sec:weakening} and detailed in Appendix.
This work has been supported by MIAI Grenoble Alpes (ANR-19-P3IA-0003).

\bibliographystyle{plainnat}
\bibliography{references}
\appendix

\section{Duality under weaker assumptions}\label{app:weakening}

In \cref{prop:general_regularized_duality}, we derive the expression of the dual of the regularized WDRO problem, with arbitrary regularizations. In this situation, the primal is not completely explicit, and then compactness and continuity are crucial, as we leverage the duality between continuous functions and measures on a compact space.

In this appendix, we provide a generalization of the duality result that does away with the compactness and the continuity assumptions. This is achieved by carefully designing the spaces in which duality arguments are made. The general duality result was informally given in \cref{sec:weakening}. Here, we provide the mathematical developments, 
as follows: \Cref{apx:setting} presents the setting; \cref{apx:state} provides the theorem; 
\cref{apx:basic} states some basic results used in the the proof of \cref{apx:proofth}; \cref{apx:examples} finishes with two examples. 

\subsection{Setting: assumptions, recalls, and construction of the spaces}\label{apx:setting}

In the appendix, we consider the following set of blanket assumptions.

\begin{assumption}%
\label{ass3}\hfill
    \begin{enumerate}[(i)]
        \item $\samples \subset \R^\dims$ is a closed set. 
        \item $\obj \from \samples \to \R$  is measurable \wrt the Borel $\sigma$-algebra on $\samples$;
        \item $\cost : \samples \times \samples \to \R_+$ is measurable as well and, for all $\sample$ in $\samples$, $\cost(\sample, \sample) = 0$;
        \item There is some $\base[\sample] \in \samples$ such that the integrability condition $\ex_{\sample \sim \prob}[\cost(\sample, \base[\sample])] < + \infty$ holds.
    \end{enumerate}
\end{assumption}

For convenience, we denote by $\base[\cost]$ the function $(\sample, \samplealt) \mapsto \cost(\sample, \base[\sample])$. 
We work in this appendix with the probability space $(\samples^2, \prodsigmaalgebra\!\!,\base)$ where  $\base$ %
is a %
measure %
satisfying the following assumption.
\begin{assumption}%
\label{ass4} $\base \in \couplings_\prob$ is a reference measure such that
        $1+\base[\cost]+\cost$ is integrable \wrt $\base$.
\end{assumption}
As a consequence of (ii), we consider $\measfuncs$ the set of measurable functions on $\samples^2$ (where functions which agree almost everywhere \wrt $\base$ are identified), equipped with the topology which metrizes the convergence in probability, see \eg \citet[Chap.~2, \S2]{kalton_peck_roberts_1984} or \citet[4.7.60]{bogachev_2007}.
For simplicity, we identify $\measfuncs$, and more generally the spaces $\lebsp[p]$, as subsets of $\measures$ through the map $\func \mapsto \func\, \dd \base$.
In turn, $\measures$ is itself seen as a subspace of $(\lebspb)^*$, the topological dual of the space of bounded (everywhere) measurable functions.

\subsubsection{Recalls in modular spaces}

The upcoming results are built over modular spaces \cite{musielak2006orlicz}, \ie subsets of $\lebsp[1]$ defined for some $\modular \from \measfuncs \to \extR$ as
\begin{align}
    \modspace \defeq  \left\{ \func\in \lebsp[1] : \exists \beta>0 , \modular(\func/\beta) < +\infty  \right\} .
\end{align}
Under suitable assumptions, the following lemma states that $\modspace$ equipped with the %
norm $ \|\func\|_\modular \defeq \inf \{ \beta>0 :  \modular(\func/\beta) \leq 1 \}$ is a Banach space.

\begin{lemma}\label{lemma:complete-modular-space}
    Consider $\modular : \measfuncs \to \R_+ \cup \{+\infty\}$ convex, symmetric, with $\modular(0) = 0$.
    If $\modular$ is $\measfuncs$-\ac{lsc} and $\setdef{\func \in \measfuncs}{\modular(\func) \leq  1}$ is bounded in $\measfuncs$, then $(\modspace, \norm{\cdot}_\modular)$ is a Banach space.
\end{lemma}
\begin{proof}
     The proof follows the second paragraph of the proof of \citet[Thm.~22]{pennanen2022topological}, which consists in applying \citet[Rmk.~9]{pennanen2022topological}. For this, one needs to show that the topology on $\modspace$ is no weaker than the $\measfuncs$ topology.
     First note that the $\measfuncs$-\ac{lsc} of $\modular$ implies the $\measfuncs$-\ac{lsc} of $\norm{\cdot}_{\modspace}$. Indeed, take a sequence $\func_\run$ of $\modspace$ which converges to $\func$:  for any $\alpha > \norm{\func_\run}_{\modspace}$, it holds that $\modular(\func_\run/\alpha) \leq 1$.

     Now, take $\open \subset \modspace$ an $\measfuncs$-neighborhood of the origin. By the boundedness assumption, there is $\alpha > 0$ such that $\setdef{\func \in \modspace\!}{\!\norm{\func}_\modular \leq  1}\!=\!\setdef{\func \in \measfuncs\!}{\!\modular(\func)\leq\!1}\!\subset\!\alpha \,\open$ so that $\open \cap \modspace$ is a neighborhood of $0$ in $\modspace$\!.\;Applying\;\citet[Rmk.\,9]{pennanen2022topological}
     yields\;the\;result.
\end{proof}

We also consider a special subspace of $\modspace$\!, named the Orlicz heart, denoted by $\heartmodspace$ and defined~by
\begin{align}
    \heartmodspace 
\defeq \left\{ \func\in \modspace : \forall \beta>0 , \modular(\func/\beta) < +\infty  \right\} .
\end{align}

Our motivation for introducing such spaces was guided by \citet{agrawal2021optimal}, which explores the link between $\phi$-divergences and associated Orlicz spaces. Note though that our construction is more general and slightly different than theirs.

\subsubsection{Construction of the function space for duality}

We first introduce two modular spaces, $\modspace[\cost]$ and $\modspace[\Reg]$, used later %
to build our function space $\pspace$.

\begin{lemma}\label{lemma:cost-modular-space}
    Define $\modular_\cost : \measfuncs \to \extR$ as
    \begin{align}
        \modular_\cost \defeq \indicator{\setdef*{\func \in \measfuncs}{ \norm*{\frac{\func}{1 + \base[\cost] + \cost}}_{\lebsp[\infty]} \leq  1 }}\,.
    \end{align}
    Then $(\modspace[\cost], \norm{\cdot}_{\modular_\cost})$ is a Banach space.
\end{lemma}
\begin{proof}
    We check that $\modular_{\cost}$ satisfies the assumptions of \cref{lemma:complete-modular-space}, and in particular that it is $\measfuncs$-\ac{lsc} and that $\setdef*{\func \in \measfuncs}{\modular_{\cost}(\func) \leq 1 }$ is bounded in probability, which comes down to showing that $ \setdef*{\func \in \measfuncs}{\norm*{{\func}/\parens*{1 + \base[\cost] + \cost}}_{\lebsp[\infty]} \leq  1 }$ is closed in $\measfuncs$ and bounded in $\measfuncs$.

    To show that it is closed, consider a sequence of functions $\func_\run$ for $\run = \running$ belonging to this set which converges in probability to some $\func \in \measfuncs$. Then, for any $\eta > 0$,
    \begin{align}
        \base \parens*{\frac{\abs{\func}}{1 + \base[\cost] + \cost} > 1 + \eta}
        &\leq 
        \base \parens*{\frac{\abs{\func - \func_\run}}{1 + \base[\cost] + \cost} > \eta}
        \base \parens*{\frac{\abs{\func_\run}}{1 + \base[\cost] + \cost} > 1} \\
        &\leq 
        \base \parens*{\abs{\func - \func_\run} > \eta} \xrightarrow[]{\run \to +\infty} 0\,.
    \end{align}
    Hence, by continuity of probability,
    \begin{align}
        \base \parens*{\frac{\abs{\func}}{1 + \base[\cost] + \cost} > 1} = \lim_{\eta \to 0}
        \base \parens*{\frac{\abs{\func}}{1 + \base[\cost] + \cost} > 1 + \eta} = 0\,,
    \end{align}
    which shows that $\norm*{{\func}/\parens*{1 + \base[\cost] + \cost}}_{L^\infty(\samples)} \leq  1$.

    We now show that $ \setdef*{\func \in \measfuncs}{\norm*{{\func}/\parens*{1 + \base[\cost] + \cost}}_{L^\infty(\samples)} \leq  1 }$ is bounded in $\measfuncs$ \ie         \begin{align}
            \forall \eta > 0,\, \exists M > 0,\, \forall \func \in \measfuncs \quad\text{ \st }\quad \, {\modular(\func) \leq 1 },\, \base(\abs{\func} \geq M) \leq  \eta\,.
        \end{align}
    To do so, let us take $\func \in \measfuncs$ such that $\norm*{{\func}/\parens*{1 + \base[\cost] + \cost}}_{\lebsp[\infty]} \leq  1$. 
    Then, for any $M > 0$, by Markov's inequality
    \begin{align}
        \base(\abs{\func} \geq M) &\leq \frac{\ex_{\base}[\abs{\func}]}{M}
                                  \leq \frac{\ex_{\base}[1 + \base[\cost] + \cost]}{M}\,,
    \end{align}
    and by assumption, $\ex_{\base}[1 + \base[\cost] + \cost]$ is finite, so that the \ac{RHS} is %
    arbitrarily small as $M$ grows.
\end{proof}

Now, let us define a modular space, starting from a ``regularization'' $\Reg$ (a more formal link will be made in \cref{apx:examples}). 

\begin{lemma}\label{lemma:reg-modular-space}
    Consider $\Reg : (\lebspb)^* \to \R_{+}$ convex such that
    \begin{enumerate}
        \item $\dom \Reg$ is included in the cone of non-negative linear forms;
        \item  $\inf \Reg = 0 = \Reg(\coupling_\Reg)$ for some $\coupling_\Reg \in\lebsp[\infty]$ such that both $(1+\alpha) \coupling_\Reg$ and $\coupling_\Reg + \alpha$ are in $\dom \Reg$ for some $\alpha > 0$;
    \end{enumerate}
    and define $\modular_\Reg : \measfuncs\to \extR$ by
    \begin{align}
        \modular_\Reg(\func) \defeq \sup_{\coupling \in \measfuncs} \int_{\samples^2} \abs{\func} \coupling \dd \base - \Reg(\coupling + \coupling_\Reg)\,.
    \end{align}
       Then $(\modspace[\Reg], \norm{\cdot}_{\modular_\Reg})$ is a Banach space. As a consequence,
       $(\heartmodspace[\Reg], \norm{\cdot}_{\modular_\Reg})$ is a Banach space too.
\end{lemma}

\begin{proof}
Note first that $\modular_\Reg(\func) = \sup_{\coupling \in \lebsp[0] : \coupling \geq 0} \int_{\samples^2}\abs{\func} \coupling \dd \base - \int_{\samples^2} \abs{\func} \coupling_\Reg \dd \base - \Reg(\coupling)$ so that it is well-defined for any $\func \in \lebsp[1]$ since the integral $\int_{\samples^2}\abs{\func} \coupling \dd \base$ always has a value in $\extR$ and because $\func \coupling_\Reg$ is still in $\lebsp[1]$.
Moreover $\modular_\Reg$ is convex as a supremum of convex functions.

Let us show that $\modular_\Reg$ satisfies the conditions of \cref{lemma:complete-modular-space}.
First, $\modular_\Reg$ is non-negative since
    \begin{align}
        \modular_\Reg(\func) \geq \inner{0}{\abs{\func}} - \Reg(0 + \coupling_\Reg) = 0\,\qquad\text{ for $\func \in \measfuncs$}.
    \end{align}
    Then, $\modular_\Reg$ is $\measfuncs$-\ac{lsc} since the functions $\func \mapsto \int{\abs{\func}}(\coupling - \coupling_\Reg)\dd \base - \Reg(\coupling)$ are $\measfuncs$-\ac{lsc} for any $\coupling \in \lebsp[0]$ non-negative by Fatou's lemma. %
    Finally, it remains to check the boundedness condition, \ie that
    $\setdef*{\func \in \measfuncs}{
                \modular_\Reg(\func)\leq 1
            }
    $
    is bounded in $\measfuncs$.
    But, since $\coupling_\Reg + \alpha \in \dom \Reg$ for some $\alpha > 0$, for any $\func \in \measfuncs$, $\modular_\Reg(\func) \geq \alpha \norm{\func}_{\lebsp[1]}$ so that $\setdef*{\func \in \measfuncs}{
                \modular_\Reg(\func)\leq 1
            }
$ is included in a $\lebsp[1]$-ball and is \emph{a fortiori} bounded in probability by Markov inequality.

 For $\heartmodspace[\Reg]$ to be Banach space too, it suffices to check that it is closed in $\modspace[\Reg]$, following the arguments of the proof of \citet[Thm.~22]{pennanen2022topological}.
 Indeed, take a sequence $(\func_\run)_\run$ from $\heartmodspace[\Reg]$ that converges to $\func \in \modspace$ for the norm $\norm{\cdot}_{\modspace[\Reg]}$. For any $\beta > 0$, $\modular_\Reg((\func - \func_\run)/\beta) \leq  1$ for $\run$ large enough so the convexity of $\modular_\Reg$ gives us that,
 \begin{align}
    \modular_\Reg(\func/(2\beta)) \leq \half \modular_\Reg(\func_\run/\beta) + \half < + \infty\,,
\end{align}
 so that $\func \in \heartmodspace[\Reg]$, and this ends the proof. 
\end{proof}

As a consequence of these two lemmas, we can build the Banach space $\pspace$ that will be used in the duality proofs, see \citet[p.\;ix]{lunardi2009interpolation} for instance.
\begin{corollary}\label{lemma:primal-space}
    Define
    $%
        \pspace = \modspace[\cost] + \heartmodspace[\Reg]
    $ %
    with the norm
    \begin{align}
        \norm{\func}_{\pspace} = \inf \setdef*{\norm{\func_\cost}_{\modspace[\cost]} + \norm{\func_\Reg}_{\modspace[\Reg]}}{\func_\cost \in \modspace[\cost],\, \func_\Reg \in \heartmodspace[\Reg]\,\st\,\, \func_\cost + \func_\Reg= \func}\,.
    \end{align}
    Then $(\pspace, \norm{\cdot}_{\pspace})$ is a Banach space, whose dual can be identified with
    \begin{align}
        \pspace^* = (\modspace[\cost])^* \cap (\heartmodspace[\Reg])^*\,.
    \end{align}
\end{corollary}

Moreover, we will need the following criterion under which a given measure belongs to $\pspace$.

\begin{lemma}\label{lemma:criterion-dual-space}
    Under the assumptions of \cref{lemma:reg-modular-space},
    for $\coupling \in \couplings_\prob \cap \measfuncs$, if both $\ex_{\coupling}[\cost]$ and $\Reg(\coupling)$ are finite, then $\coupling$ belongs to $\dspace$.
\end{lemma}
\begin{proof}
    We first show that $\coupling$ is in $(\modspace[\cost])^*$. To do so, we use that by \citet[Thm.~22]{pennanen2022topological}, $\kdmodspace[\cost]$ is a subspace of $(\modspace[\cost])^*$ where $\modular_\cost^*$ is defined for all $\coupling\in(\modspace[\cost])^*$ as  $\modular_\cost^*(\coupling) = \sup_{\func \in \lebsp[\infty]} \inner{\coupling}{\func} - \modular_\cost(\func) $. 
    By definition of $\modular_\cost$, we have that $\modular_\cost^*(\coupling) = \inner{\coupling}{1 + \base[\cost] + \cost } = 1 + \ex_{\sample \sim \prob}[\cost(\sample, \base[\sample])] + \ex_{\coupling}[\cost]$, which is finite by assumption. We thus have that $\coupling \in (\modspace[\cost])^*$. 

    Next, we show that $\coupling$ defines a continuous linear form on $(\heartmodspace[\Reg], \norm{\cdot}_{\modspace[\Reg]})$. To do so, let us take $\func \in \heartmodspace[\Reg]$ so that by definition we can take $\beta > 0$ so that $\modular_\Reg(\func/\beta) < +\infty$. 
    By the assumption of \cref{lemma:reg-modular-space}, there is some $\alpha > 0$ such that $(1 + \alpha) \coupling_\reg \in \dom \Reg$ so that one can choose $\lambda \in (0, 1)$ satisfying  $(1+\alpha)(1-\lambda) = 1$.
    Then, by definition of $\modular_\Reg$ and convexity of $\Reg$, one has that for all $\coupling\in  \measfuncs$, the following inequality holds
    \begin{align}\label{eq:app-proof-criterion-dual}
        \frac{\lambda}{\beta} \int \abs{\func} \coupling \dd \base \leq \modular_\Reg(\func / \beta) + \lambda \Reg(\coupling) + (1 - \lambda) \Reg((1 + \alpha) \coupling_\Reg)\,.
    \end{align}
    Since the \ac{RHS} is finite from the assumptions as soon as $\Reg(\coupling)$ is finite, we have that $\int {\func} \coupling \dd \base$ is well-defined for $\func \in \heartmodspace[\Reg]$. All that is left to show is that this linear form is indeed continuous. Take $\func_\run \in \heartmodspace[\Reg]$ which goes to zero as $\run$ goes to $+\infty$. This means that there exists a sequence $\beta_\run> 0$ which goes to zero such that $\modular_\Reg(\func_\run/\beta_\run) \leq 1$ for all $\run = \running$.
    Applying \cref{eq:app-proof-criterion-dual} with $\func \gets \func_\run$, $\beta \gets \beta_\run$, and $\coupling$ such that $\Reg(\coupling)$ is finite yields
   \begin{align}
        \frac{\lambda}{\beta_\run} \int \abs{\func_\run} \coupling \dd \base \leq 1 + \lambda \Reg(\coupling) + (1 - \lambda) \Reg((1 + \alpha) \coupling_\Reg) < +\infty    \end{align}
    which implies that $\int \abs{\func_\run} \coupling \dd \base $ goes to 0 as $\run$ goes to $+\infty$.
\end{proof}

\subsection{Statement of the theorem}\label{apx:state}

We can now state the formal version of \cref{theorem:weakened_general_regularized_duality}. 

\begin{theorem}\label{theorem:app_weakened_general_regularized_duality}
    Take $\Regalt \defeq \frac{\regalt}{\reg} \Reg$.
    Let \cref{ass3} hold and assume that
    \begin{enumerate}
        \item $\obj$ satisfies the growth condition :
            \begin{align}
                \exists \func \in \lebsp[\infty],\,
                \forall \beta > 0,\,
                (\sample, \samplealt) \mapsto \beta \left(\obj(\samplealt) - (1 + \cost(\sample, \base[\sample]) + \cost(\sample, \samplealt)) \func(\sample, \samplealt)\right) \in \dom \modular_\Reg\,,
            \end{align}
            or equivalently, that the function $\tildeobj$, defined as $\tildeobj : (\sample, \samplealt) \mapsto \obj(\samplealt)$, lies in $\pspace$;
        \item $\Reg$ convex satisfies the assumptions of \cref{lemma:reg-modular-space} and $\Reg \restriction_{\dspace}$ is $\weaktop{\dspace}{\pspace}$-\ac{lsc};
        \item $\base \in \couplings_\prob$ satisfies \Cref{ass4};
        \item there exists 
$\coupling_\Regalt \in \couplings_\prob \cap \measfuncs$ such that $\ex_{\coupling_\Regalt}[\cost] + \Regalt(\coupling_\Regalt) < \radius$; %
        \item one of the two following regularity conditions hold:
 \begin{enumerate}
        \item $\Reg$ is $\weaktop{\dspace}{\contfuncs \cap \modspace[\cost]}$-\ac{usc} on $\couplings_\prob$, $\obj$ is \ac{lsc}, $\obj_- \in \modspace[\cost]$, $\overline{\intr \samples} = \samples$ and for all $\sample \in \samples$ and $\base(\cdot | \sample)$ has a positive density \wrt the Lebesgue measure on $\samples$.
        \item $\dom \Reg \subset \measfuncs$ and, for any non-increasing non-negative sequence  $\func_\run \in \dom \modular_\Reg$  such that $\func_\run \to 0$ as $\run \to \infty$, we have $\modular_\Reg(\func_\run) \to 0$.
    \end{enumerate}
    \end{enumerate}
    	Then, the problem
	\begin{align}
		\sup \left\{
		\ex_{\coupling_2}[\obj] - \Reg(\coupling)
		:
		\coupling \in \couplings_\prob\,, \ex_{\coupling}[\cost] + \Regalt(\coupling) \leq \radius
		\right\}
	\end{align}
	coincides with its dual
	\begin{align}
		\inf \left\{
        \dualvar \radius +
		\ex_{\sample \sim \prob} \left[
            \esssup
        \left\{
			\obj - \dualvar \cost(\sample, \cdot) - \func(\sample, \cdot)
        \right\}
        \right]
        +
        \left(\Reg\restriction_{\dspace} + \dualvar \Regalt\restriction_{\dspace}\right)_*(\func)
        :
		\func \in \pspace,\, \dualvar \geq 0
		\right\}\,.
	\end{align}
\end{theorem}

We state the theorem and provide its proof in the case of when $\Reg$ and $\Regalt$ are proportional. The case with different $\Reg$ and $\Regalt$ would follow that the same lines, at the price of slightly more complicated assumptions and a slightly different treatment of subcases in the proof. To avoid these extra-technicalities, we stick here with the proportional case, which already presents all the intrinsic reasoning and the mathematical tools. 
Before going to the proof of this result in \cref{apx:proofth}, we introduce in the next section some general duality lemmas, tailored for our context.

\subsection{Basic duality lemmas}\label{apx:basic}

\begin{lemma}\label{lemma:biconjugate}
	Consider $\pspace$ a Banach space and $\func \from \dspace \to \extR$ convex, proper and weakly-$\star$ \ac{lsc}. Then $(\func_*)^* = \func$.
\end{lemma}
\begin{proof}
	Equipped with its weak-$\star$ topology, $\dspace$ is a locally convex topological vector space whose topological dual is $\pspace$, see \citet[Thm.~3.10, \S 3.14]{rudin1991functional}. As a consequence, the usual theorem of the biconjugate (\eg \citet[Thm.~2.3.3]{zalinescu2002convex}) ensures that $(\func_*)^* = \func$.
\end{proof}

\begin{lemma}\label{lemma:duality}
	Given $\pspace$ a Banach space, $\func, \funcalt, \constr, \constralt \from \dspace \to \extR$ convex, proper and weakly-$\star$ \ac{lsc}, consider the two problems,
	\begin{align}
		\min_{\point \in \pspace,\, \dualvar \geq 0} (\func + \dualvar \constr)_*(\point) + (\funcalt + \dualvar \constralt)_*(-\point) \tag{P}\label{eq:primal_abstract} \\
		\max_{\dpoint \in \dspace} -\func(\dpoint) -\funcalt(\dpoint) \text{ s.t. } \constr(\dpoint)+\constralt(\dpoint) \leq  0 \tag{D}\label{eq:dual_abstract}\,.
	\end{align}
	If the following qualification condition holds
    \begin{align}
   \setdef*{\point \in \pspace}{\exists \dualvar \geq 0,\, \point \in \dom (\func + \dualvar \constr)_* + \dom (\funcalt + \dualvar \constralt)_*} = \pspace
    \end{align}
	then, $\cref{eq:primal_abstract} = \cref{eq:dual_abstract}$.
\end{lemma}
\begin{proof}
	Let us define the perturbation function $\pert : (\pspace \times \R) \times \pspace \to \extR$ by
	\begin{align}
		\pert((\point, \dualvar), \pointalt) = (\func + \dualvar \constr)_*(\point) + (\funcalt + \dualvar \constralt)_*(\pointalt - \point) + \indicator{\R_+}(\dualvar)\,.
	\end{align}
	$\pert$ is \ac{lsc} and convex since both $(\point, \dualvar), \pointalt \mapsto (\func + \dualvar \constr)_*(\point)$ and $(\point, \dualvar), \pointalt \mapsto (\funcalt + \dualvar \constralt)_*(\pointalt - \point)$ are jointly convex and \ac{lsc} as supremums of continuous linear functions.
	Let us now compute the conjugate of $\pert$, for $((\dpoint, \ddualvar), \dpointalt) \in (\dspace \times \R) \times \dspace$
	\begin{align}
		\pert^*((\dpoint, \ddualvar), \dpointalt)
		 & = \sup_{\point, \pointalt \in \pspace,\, \dualvar \geq 0} 
         \left\{
         \inner{\dpoint, \point} + \inner{\dpointalt, \pointalt} + \ddualvar \dualvar - (\func + \dualvar \constr)_*(\point) - (\funcalt + \dualvar \constralt)_*(\pointalt-\point)  
     \right\} \\
		 & = \sup_{\point, \pointalt \in \pspace,\, \dualvar \geq 0} 
        \left\{
         \inner{\dpoint, \point} + \inner{\dpointalt, \pointalt + \point} + \ddualvar \dualvar - (\func + \dualvar \constr)_*(\point) - (\funcalt + \dualvar \constralt)_*(\pointalt) 
     \right\}\\
		 & = \sup_{\dualvar \geq 0}  
         \left\{
         \ddualvar \dualvar + (\func + \dualvar \constr)(\dpoint + \dpointalt) + \sup_{\point \in \pspace} \{ \inner{\dpoint, \point}-  (\funcalt + \dualvar \constralt)_*(\pointalt) \} 
     \right\} \\
		 & = \sup_{\dualvar \geq 0} \left\{
             \ddualvar \dualvar + (\func + \dualvar \constr)(\dpoint + \dpointalt) +  (\funcalt + \dualvar \constralt)(\dpointalt)                                                                                      
             \right\} \\
		 & = \func(\dpoint + \dpointalt) + \funcalt(\dpointalt) + \indicator{\ddualvar + \constr(\dpoint + \dpointalt) + \constralt(\dpointalt) \leq 0}\,,
	\end{align}
	where we applied \cref{lemma:biconjugate} to $\func + \dualvar \constr$ and $\funcalt + \dualvar \constralt$.
	Now, note that \cref{eq:primal_abstract} is equal to
	\begin{align}
		\min_{(\point, \dualvar) \in (\pspace, \R)} \pert((\point, \dualvar), 0)
        \qquad\text{and \cref{eq:dual_abstract} to}\qquad \max_{\dpoint \in \dspace} - \pert^*(0, \dpoint)\,.
	\end{align}
	Now we leverage the duality theorem \citet[Thm.~2.7.1]{zalinescu2002convex} with qualification condition (vii), see \citet[p.~15]{zalinescu2002convex} for the relevant definition, to show that the values of those problems are equal. This qualification condition requires that the set
	\begin{align}
		\cone \{\pointalt \in \pspace : \exists (\point, \dualvar) \in \pspace \times \R,\, ((\point, \dualvar), \pointalt) \in \dom \pert \}
	\end{align}
	be a closed linear subset of $\pspace$. But this set can be rewritten as
	\begin{align}
		\cone \{\pointalt \in \pspace : \exists \dualvar \geq 0,\, \pointalt \in \dom (\func + \dualvar \constr)_* + \dom (\funcalt + \dualvar \constralt)_* \}
	\end{align}
	and our qualification assumption ensures that it is equal to the whole space $\pspace$.
\end{proof}

\begin{lemma}\label{lemma:ineq-constraints}
    Consider a vector space $\dspacealt$, $\func : \dspacealt \to \{-\infty\}\cup \R$ concave, $\constr : \dspacealt \to \extR$ convex and assume that there exists $\dpointalt_\constr \in \dom \func$ such that $\constr(\dpointalt_\constr) < 0$. Then,
    \begin{align}
        \sup \setdef{\func(\dpointalt)}{\dpointalt \in \dspacealt,\, \constr(\dpointalt) \leq  0}
        =
        \sup \setdef{\func(\dpointalt)}{\dpointalt \in \dspacealt,\, \constr(\dpointalt) < 0}\,.
    \end{align}
\end{lemma}
\begin{proof}
    The inequality $(\geq)$ always holds. To prove the reverse, take $\dpointalt \in \dom \func$ such that $\constr(\dpointalt) \leq 0$, 
    consider 
    $\dpointalt_t \defeq t \dpointalt_\constr + (1-t) \dpointalt$ for $t \in [0, 1]$. By convexity of $\constr$ and definition of $\dpointalt_\constr$, we have that $\constr(\dpointalt_t) < 0$ for all $t \in (0,1]$. Moreover, since $\dpointalt_\constr$ is in $\dom \func$ and $\func$  is concave, it holds that %
    \begin{align}
        \func(\dpointalt) 
        \leq \liminf_{t \to 0} \frac{\func(\dpointalt_t) - t \func(\dpointalt_\constr)}{1 - t}
        = \liminf_{t \to 0} \func(\dpointalt_t)
        \leq \sup \setdef{\func(\dpointalt)}{\dpointalt \in \dspacealt,\, \constr(\dpointalt) < 0}\,,
    \end{align}  
    which shows the required inequality.
\end{proof}

\subsection{Proof of \cref{theorem:app_weakened_general_regularized_duality}}
\label{apx:proofth}

The proof of \cref{theorem:app_weakened_general_regularized_duality} is divided two main lemmas, \cref{lemma:proof_app_partone} and \cref{lemma:reformulate-dual}, which, put together, exactly give \cref{theorem:app_weakened_general_regularized_duality}.
\begin{lemma}\label{lemma:support-function}
	Assume that $\pspace \subset \lebsp[1]$ is a Banach space 
	then the support function of the set $\couplings_\prob \cap \measfuncs \cap \dspace$, seen as a subspace of $\dspace$, is given for all  $\func \in \pspace$ by
	\begin{align}
		\support_{\couplings_\prob \cap \measfuncs \cap \dspace}(\func)  \defeq
		\sup_{\coupling \in \couplings_\prob \cap \measfuncs \cap \dspace} \inner{\coupling, \func} 
		                                                     =
		\ex_{\sample \sim \prob} [
			\esssup \func(\sample, \cdot)
		]
	\end{align}
\end{lemma}

\begin{remark}\label{remark:ess-sup}
	We observe  that
	$
		\ex_{\sample \sim \prob} [
			\esssup \func(\sample, \cdot)
		]
	$ is well-defined for $\func \in \pspace$. %
        First, the map $\sample \mapsto \esssup \func(\sample, \cdot)$ is measurable since, for any $\const \in \R$, by definition of the essential supremum,
        \begin{align}
            \setdef*{\sample \in \samples}{\esssup \func(\sample, \cdot) > \const}
            &=
            \setdef*{\sample \in \samples}{ \base \parens*{ \setdef{\samplealt \in \samples}{\func(\sample, \samplealt) > \const} | \sample} > 0}\\
            &=
            \setdef*{\sample \in \samples}{ \int_{\samples} \one_{\setdef{(\sample, \samplealt) \in \samples^2}{\func(\sample,\samplealt) > \const}}(\sample, \samplealt)\dd \base(\samplealt | \sample) > 0}\,,
        \end{align}
        where the resulting integral is a measurable function of $\sample$, see \eg \citet[Lem.~3.2 (i)]{kallenberg2017random}.
    Then, if $\func \in \pspace$, then for $\base$-almost every $(\sample, \samplealt)$,
			$\esssup \func(\sample, \cdot) \geq \func(\sample, \samplealt)$
        so that by integrating \wrt $\base$, one gets that 
    $
        \ex_{\sample \sim \prob} [
			\esssup \func(\sample, \cdot)
		]
        \geq \ex_{\base}[\func] > - \infty
    $
    since $\func \in \pspace \subset \lebsp[1]$.
\end{remark}

\begin{proof}
Fix some $\func \in \pspace$. 
First, let us show the inequality $(\leq)$. To do so, take $\coupling \in \couplings_\prob \cap \measfuncs \cap \dspace$. By disintegration of measure, it holds that
\begin{align}
    \ex_{\coupling}[\func]
    = \ex_{\sample \sim \prob}[\ex_{\samplealt \sim \coupling(\cdot | \sample)}[\func(\sample, \samplealt)]] \leq  \ex_{\sample \sim \prob}[\esssup{\func(\sample, \cdot)}]\,,
\end{align}
since $\coupling(\cdot|\sample)$ is a probability distribution which is absolutely continuous \wrt $\base[\coupling](\cdot | \sample)$ for any $\sample \in \samples$.

To prove the reverse inequality, we distinguish three cases. For convenience, denote by $\funcalt$ the measurable function $\sample \mapsto \esssup{\func(\sample, \cdot)}$. 
Furthermore, note that the fact that $\pspace$ is a subset of $\lebsp[1]$ implies that $\dspace$ contains $\lebsp[\infty]$, so the indicator functions of measurable sets, in particular.
\begin{enumerate}
\item If $\parens*{\prob-\als ~ \func(\sample, \cdot) \in L^\infty(\base(\cdot|\sample))}$ is false, \ie if $\funcalt$ is not finite $\prob-\als$.
    This means that, for any $n \geq 1$, the sets $\mset_n(\sample) \defeq \setdef{\samplealt \in \samples}{\func(\sample, \samplealt)\geq n}$ satisfy $\prob-\als,\, \base(\mset_n(\sample)|\sample) > 0$.
    Therefore, one can define $\coupling_n \in \couplings_\prob \cap \measfuncs \cap \dspace$ as the measurable map $(\sample, \samplealt) \mapsto \frac{\one_{\mset_n(\sample)}(\samplealt)}{\base(\mset_n(\sample)|\sample)}$ and it satisfies $\inner{\coupling_n}{\func} \geq n$.
    This shows that
    \begin{align}
	\sup_{\coupling \in \couplings_\prob \cap \measfuncs \cap \dspace} \inner{\coupling, \func} = +\infty
	=
		\ex_{\sample \sim \prob} [
                \funcalt(\sample)
		]\,.
    \end{align}

\item If ${\prob-\als ~  \func(\sample, \cdot) \in L^\infty(\base(\cdot|\sample))}$, \ie if  $\funcalt$ is finite $\prob-\als$ \emph{and} $\funcalt \in L^1(\prob)$.
 As a consequence of the definition of $\funcalt$, for any $n \geq 1$, the sets $\mset_n(\sample) \defeq \setdef{\samplealt \in \samples}{\func(\sample, \samplealt)\geq \funcalt(\sample) - 1/n}$ satisfy $\prob\als,\, \base(\mset_n(\sample)|\sample) > 0$.
    Again, one defines $\coupling_n \in \couplings_\prob \cap \measfuncs \cap \dspace$ as the %
    map $(\sample, \samplealt) \mapsto \frac{\one_{\mset_n(\sample)}(\samplealt)}{\base(\mset_n(\sample)|\sample)}$ so that it satisfies $\inner{\coupling_n}{\func} \geq \ex_{\prob}[\funcalt]- 1/n$ for all $n \geq 1$ and thus
    \begin{align}
	\sup_{\coupling \in \couplings_\prob \cap \measfuncs \cap \dspace} \inner{\coupling, \func}
	\geq
		\ex_{\sample \sim \prob} [
                \funcalt(\sample)
		]\,.
    \end{align}
\item Finally, if ${\prob-\als ~  \func(\sample, \cdot) \in L^\infty(\base(\cdot|\sample))}$, \ie if  $\funcalt$ is finite $\prob-\als$ \emph{but} $\funcalt \notin L^1(\prob)$.
    By \cref{remark:ess-sup}, this means that $\ex_{\prob}[\funcalt] = +\infty$.
    By mimicking the proof of the previous point but with $\funcalt_n(\sample, \samplealt) \defeq \funcalt(\sample) \one_{\funcalt(\sample) \leq n} + \func(\sample, \samplealt) \one_{\funcalt(\sample) > n}$ instead of $\funcalt$, we obtain that, for any $n \geq 1$,
\begin{align}
    \sup_{\coupling \in \couplings_\prob \cap \measfuncs \cap \dspace} \inner{\coupling, \func}
	\geq
		\ex_{\sample \sim \base} [
                \funcalt_n(\sample, \samplealt)
		]\,. %
\end{align}
    Now, by applying the monotone convergence theorem to $\funcalt_n - \func$, which is non-negative $\base-\als$, and with $\funcalt \in \lebsp[1]$  by assumption, we get that
\begin{align}
    \sup_{\coupling \in \couplings_\prob \cap \measfuncs \cap \dspace} \inner{\coupling, \func}
	\geq
		\ex_{\sample \sim \prob} [
                \funcalt(\sample)
		] = +\infty\,, %
\end{align}
which concludes the proof.
\end{enumerate}
\vspace*{-3ex}
\end{proof}

\begin{lemma}\label{lemma:proof_app_partone}
    Under the assumptions of \cref{theorem:app_weakened_general_regularized_duality},
	the problem
	\begin{align}
		\sup \left\{ \ex_{\coupling_2}[\obj]
	 - \Reg(\coupling)
		:
		\coupling \in \wstarcl{\couplings_\prob \cap \measfuncs\cap\dspace}\,, \inner{\coupling}{\cost} + \Regalt(\coupling) \leq \radius
		\right\}
	\end{align}
	coincides with its dual
	\begin{align}
		\inf \left\{
        \dualvar \radius +
		\ex_{\sample \sim \prob} \left[
            \esssup
        \left\{
			\obj - \dualvar \cost(\sample, \cdot) - \func(\sample, \cdot)
        \right\}
        \right]
        +
        \left(\Reg\restriction_{\dspace} + \dualvar \Regalt\restriction_{\dspace}\right)_*(\func)
        :
		\func \in \pspace,\, \dualvar \geq 0
		\right\}\,.
	\end{align}
\end{lemma}

\begin{proof}
    We apply \cref{lemma:duality} with $\pspace$ and $\func \gets \Reg\restriction_{\dspace}$, $\constr \gets \Regalt\restriction_{\dspace}$ which are convex, proper and $\weaktop{\dspace}{ \pspace}$-\ac{lsc}. Indeed,  %
    $\weaktop{\dspace}{\modspace[\cost]}$ is weaker than $\weaktop{\dspace}{\pspace}$, and $\funcalt \gets \inner*{\cdot}{- \tildeobj} + \indicator{\wstarcl{\couplings_\prob \cap \measfuncs \cap \dspace}}$, $\constralt \gets \inner*{\cdot}{\cost} - \radius$ are convex, proper and weakly-$\star$ \ac{lsc} by construction.
    For $\dualvar \geq 0$, we have 
    \begin{align}
        \forall \func \in \pspace,\quad (\funcalt + \dualvar \constralt)_*(\func) &= \dualvar \radius +\sup_{\coupling \in \wstarcl{\couplings_\prob \cap \measfuncs\cap\dspace}} \inner*{\coupling}{\func} - \inner*{\coupling}{\dualvar \cost - \tildeobj} \\
                                                                                             &= \dualvar \radius + \support_{\wstarcl{\couplings_\prob \cap \measfuncs\cap\dspace}}\left(\tildeobj - \dualvar \cost + \func\right) \\
                                                                                             &= \dualvar \radius + \support_{{\couplings_\prob \cap \measfuncs\cap\dspace}}\left(\tildeobj - \dualvar \cost + \func\right)\\
                                            &=  \dualvar \radius +
		\ex_{\sample \sim \prob} \left[
            \esssup
        \left\{
			\obj - \dualvar \cost(\sample, \cdot) - \func(\sample, \cdot)
        \right\}
        \right]
    \end{align}
    where the last equality comes from \cref{lemma:support-function}. The qualification condition of \cref{lemma:duality} writes,
    \begin{align}\label{eq:gen-duality-cost-qc}
   \setdef*{\func \in \pspace}{\exists \dualvar \geq 0,\, \func \in \dom (\Reg + \dualvar \Regalt)_* + \dom \support_{{\couplings_\prob \cap \measfuncs\cap\dspace}} - \left(\tildeobj - \dualvar \cost\right)} = \pspace\,.
    \end{align}
Take any $\func \in \pspace$ and let us show that it actually belongs to this set. Since $\tildeobj$ also belongs to $\pspace$ by assumption, there exists $\funcalt \in \heartmodspace[\Reg]$ \and some constant $\const > 0$ such that, $\base[\coupling]$ everywhere
\begin{align}
        \func(\sample, \samplealt) + \tildeobj(\sample, \samplealt) - \funcalt(\sample, \samplealt) \leq  \const \parens*{1 + \base[\cost](\sample, \samplealt) + \cost(\sample, \samplealt)}\,.
\end{align}
Hence, with $\dualvar \defeq \const$,
\begin{align}
        \func(\sample, \samplealt) + \tildeobj(\sample, \samplealt) - \funcalt(\sample, \samplealt)- \const\cost(\sample, \samplealt) \leq  \const \parens*{1 + \cost(\sample, \base[\sample])}\,.
\end{align}
But the \ac{RHS} belongs to $\dom \support_{\couplings_\prob \cap \kdmodspace}$ since $\ex_{\sample \sim \prob}[\cost(\sample, \base[\sample])] < + \infty$. Then the \ac{LHS} is inside the domain, as well. Moreover, $\funcalt$ belongs to $\heartmodspace[\Reg]$ so that it lies inside $\dom (\Reg + \dualvar \Regalt)_* = \frac{1}{1 + \dualvar \regalt / \reg}\dom \Reg$.
Hence, we have shown that $\func$ belongs to the \ac{LHS} of~\cref{eq:gen-duality-cost-qc}, which shows the equality.
\end{proof}

\begin{lemma}
\label{lemma:dense}
    Assume that $\cost$ is continuous on $\samples^2$.
    For any $\coupling \in \couplings_\prob$ such that $\ex_{\coupling}[\cost]$ is finite, there exists $\coupling_\run \in \couplings_\prob$ for $\run = \running$ such that
    \begin{enumerate}
    \item  $\coupling_\run$ is absolutely continuous \wrt $\base$ for $\run = \running$,
    \item  $(\coupling_\run)_\run$ converges weakly to $\coupling$, \ie \wrt the topology $\weaktop{\measures}{\bdedcontfuncs}$,
    \item  $\ex_{\coupling_\run}[\cost] \to \ex_{\coupling}[\cost]$ as $\run \to \infty$.
    \end{enumerate}
   In particular, for any $\func \in \measfuncs$ \ac{lsc} such that $\func_- \in \modspace[\cost]$,
   \begin{align}
       \ex_{\coupling}[\func] \leq  \liminf_{\run \to \infty} \ex_{\coupling_\run}[\func]\,,
   \end{align}
   and, as a consequence, $(\coupling)_\run$ also converges to $\coupling$ for the topology $\weaktop{\dmodspace[\cost]}{\contfuncs \cap \modspace[\cost]}$.
\end{lemma}
\begin{proof}

\noindent\underline{Step 1: building a partition of $\samples^2$.}
For convenience, denote by $\norm{\cdot}$ a norm on $\R^{\dims}$, which contains $\samples$ and by $\ball(0, r)$ the closed ball of radius $r$ around 0 in $\R^\dims$.
Fix $\Regparam > 0$. By absolute continuity of $\cost$ on compact sets, for any $n, m \geq 1$, there exists $\eta_{n, m} \in (0, \Delta]$ such that, for any $\sample \in \samples \cap \ball(0, n)$, $\samplealt,\,\samplealt' \in \samples \cap \ball(0, m)$,
\begin{align}
    \norm{\samplealt - \samplealt'} \leq \eta_{n, m} \implies \abs{\cost(\sample, \samplealt')- \cost(\sample, \samplealt)} \leq \Regparam\,.
\end{align}

Fix $n$ and $m$ greater or equal than 1.%
We consider a finite family of disjoint open sets $\open_i^m$ for $i \in I^m$ which are included in $\openball(0, m) \setminus \ball(0, m-1)$, whose diameters are at most $\eta_{n,m}$ and which satisfy $\overline{\bigcup_{i \in I^m} \open_i^m} \supset \openball(0, m) \setminus \ball(0, m-1)$.
Define $J^m = \setdef{i \in I^m}{\open_i^m \cap \intr \samples \neq \emptyset}$. We first show%

\begin{align}\label{eq:app-geometric-open-sets-partition}
    \overline{\bigcup_{i \in J^m} \open_i^m \cap \intr \samples} \supset \left( \openball(0, m) \setminus \ball(0, m-1) \right) \cap \samples\,
\end{align}
by considering $\samplealt \in (\openball(0, m) \setminus \ball(0, m-1)) \cap \samples$ and $\open$ neighborhood of $\samplealt$. We thus have to prove that $\bigcup_{i \in J^m} \open_i^m \cap \intr \samples \cap \open$ is not empty. 
Since $\openball(0, m) \setminus \ball(0, m-1)$ is an open set that contains\;$\samplealt$, it suffices to show this for some $\open \subset (\openball(0, m) \setminus \ball(0, m-1))$. This means that, in particular, $\open \subset \overline{\bigcup_{i \in I^m} \open_i^m}$. Because $\samplealt$ belongs to $\samples$ with $\samples = \overline{\intr \samples}$, it holds that $\open \cap \intr \samples \neq \emptyset$.
Therefore, we have that $\open \cap \intr \samples \cap \overline{\bigcup_{i \in I^m} \open_i^m} \neq \emptyset$. Since $\open \cap \intr \samples$ and $\bigcup_{i \in I^m} \open_i^m$ are both open, this implies that $\open \cap \intr \samples \cap {\bigcup_{i \in I^m} \open_i^m} \neq \emptyset$, which concludes the proof of \eqref{eq:app-geometric-open-sets-partition}.
Moreover, since $J^m$ is finite, 
$
\overline{\bigcup_{i \in J^m} \open_i^m \cap \intr \samples} 
={\bigcup_{i \in J^m} \overline{\open_i^m \cap \intr \samples}} 
$ so that 
$
{\bigcup_{i \in J^m} \overline{\open_i^m \cap \intr \samples}} 
\supset \left( \openball(0, m) \setminus \ball(0, m-1) \right) \cap \samples\,
$ holds as well.

Still keeping $n \geq 1$ fixed, we have built an (at most) countable number of disjoint open sets $(\openalt_i)_{i \leq M_n} \defeq(\open^m_i \cap \intr \samples)_{i \in J^m,\, m \geq 1}$, with $M_n \in \mathbb{N} \cup \{+\infty\}$ such that for any $i \leq  M_n$, there is some $m \geq 1$ such that $\openalt_i \subset \openball(0, m) \setminus \ball(0, m-1)\cap \intr \samples$ and such that the diameter of $\openalt_i$ is at most $\eta_{n, m}$. Moreover, $\bigcup_{i \in I} \overline{\openalt_i}$ is the whole of $\samples$. 
We build a partition of $\samples$ by recursively defining
\begin{align}
        \openaltalt^n_{i} = \overline{\openalt_i} \setminus \bigcup_{j \leq i-1} \openaltalt^n_j
        = \openalt_i \cup \parens*{\partial \openalt_i \setminus \bigcup_{j \leq i-1} \openaltalt^n_j} \text{ for } i \leq M_n\,.
\end{align}
Therefore, we have built an (at most) countable number of disjoint measurable sets $\openaltalt^n_i$ for $i \leq M_n$, whose union is $\samples$ and such that for any $i \leq  M_n$, there is some $m \geq 1$ such that $\openaltalt_i^n \subset \ball(0, m) \cap \samples$. Moreover, the diameter of each $\openaltalt_i^n$ is at most $\eta_{n, m}$ and each of these sets has non-empty interior.

The full partition of $\samples^2$ that we finally consider is made of the sets $\parens*{(\openball(0, n) \setminus \openball(0, n-1)) \cap \samples} \times \openaltalt^n_i$ for $n \geq 1$ and $i \leq M_n$. To summarize, we have built, for $\Delta > 0$, a partition of the form $(\mset_i^\Delta \times \msetalt_j^\Delta)_{i \in I,\, j \in J(i)}$ which is a (at most) countable and measurable partition of $\samples^2$ satisfying:%
\begin{enumerate}
    \item $(\mset_i^\Delta)_{i \in I}$ is a partition of $\samples$ and, for any $i \in I$, $(\msetalt_j^\Delta)_{j \in J(i)}$ is a partition of $\samples$.
    \item For $i \in I$, $j \in J(i)$, $\sample \in \mset_i^\Delta$, and $\samplealt, \samplealt' \in \msetalt_j^\Delta$, we have $\norm{\samplealt - \samplealt'} \leq \Delta$ and $\abs{\cost(\sample, \samplealt')- \cost(\sample, \samplealt)} \leq \Delta$.
    \item For any $i \in I$, $j \in J(i)$, $\msetalt_j^\Delta$ has non-empty interior.
\end{enumerate}

\noindent\underline{Step 2: construction of the sequence of measures.}
For $\Delta > 0$, we define the measure $\coupling^\Delta \in \couplings_\prob$ for any $ \sample \in \samples$ and any measurable $\msetalt \subset \samples $ by
\begin{align}
    \coupling^\Delta(\msetalt | \sample) &\defeq \sum_{i \in I} \one_{\sample \in \mset_i^\Delta}\sum_{j \in J(i)} \frac{\coupling(\msetalt_j^\Delta | \sample)}{\base(\msetalt_j^\Delta | \sample)} \base(\msetalt \cap \msetalt_j^\Delta | \sample)\,,
\end{align}
which is well-defined because of the $B_j^\Delta$ having non-empty interiors. Moreover, by construction, it is absolutely continuous \wrt $\base$, so that $\coupling^\Delta$ satisfies $(1)$.

For\;$i \in I$ and $j \in J(i)$, choose $\samplealt_{i,j}^\Delta \in \msetalt_j^\Delta$.
By\;the\;Portmanteau\;theorem, see\;\citet[Th\,13.16]{klenke2014probability}, it suffices to show that $\int \func \, \dd \coupling^\Delta$ converges to $\int \func \, \dd \coupling$ as $\Delta \to 0$ for any bounded Lipschitz function $\func$ to prove the weak convergence of the sequence of measures $(\coupling^\Delta)$ to $\coupling$. Thus, take a bounded $1$-Lipschitz function $\func : \samples^2 \to \R$, then 
\begin{align}
    \abs*{\int \func \, \dd \coupling^\Delta - \int \func \, \dd \coupling}
    \leq& \ex_{\sample \sim \prob}[
        \sum_{i \in I} \one_{\sample \in \mset^\Delta_i}
        \sum_{j \in J(i)} 
        \abs*{
        \frac{\coupling(\msetalt_j^\Delta | \sample)}{\base(\msetalt_j^\Delta | \sample)}
        \int_{\msetalt_j^\Delta} \func(\sample, \samplealt) \base(\dd \samplealt | \sample)
        - \int_{\msetalt_j^\Delta} \func(\sample, \samplealt) \coupling(\dd \samplealt | \sample)
    }
    ]\\
    \leq& \ex_{\sample \sim \prob}[
        \sum_{i \in I} \one_{\sample \in \mset^\Delta_i}
        \sum_{j \in J(i)} 
        \coupling(\msetalt_j^\Delta | \sample)
        \parens*{
        \frac{1}{\base(\msetalt_j^\Delta | \sample)}
        \int_{\msetalt_j^\Delta} \abs*{\func(\sample, \samplealt) - \func(\sample, \samplealt_{i,j}^\Delta)} \base(\dd \samplealt | \sample)
    }]\\
    &+\ex_{\sample \sim \prob}[
        \sum_{i \in I} \one_{\sample \in \mset^\Delta_i}
        \sum_{j \in J(i)} 
        \coupling(\msetalt_j^\Delta | \sample)
        \parens*{
        \frac{1}{\coupling(\msetalt_j^\Delta | \sample)}
        \int_{\msetalt_j^\Delta} \abs*{\func(\sample, \samplealt) - \func(\sample, \samplealt_{i,j}^\Delta)} \coupling(\dd \samplealt | \sample)
    }]
    \\
    \leq& \ex_{\sample \sim \prob}[
        \sum_{i \in I} \one_{\sample \in \mset^\Delta_i}
        \sum_{j \in J(i)} 
        \coupling(\msetalt_j^\Delta | \sample)
        \parens*{
        \frac{1}{\base(\msetalt_j^\Delta | \sample)}
        \int_{\msetalt_j^\Delta} {\norm{\samplealt - \samplealt_{i,j}^\Delta}} \base(\dd \samplealt | \sample)
    }]\\
    &+\ex_{\sample \sim \prob}[
        \sum_{i \in I} \one_{\sample \in \mset^\Delta_i}
        \sum_{j \in J(i)} 
        \coupling(\msetalt_j^\Delta | \sample)
        \parens*{
        \frac{1}{\coupling(\msetalt_j^\Delta | \sample)}
        \int_{\msetalt_j^\Delta} {\norm{\samplealt - \samplealt_{i,j}^\Delta}}  \coupling(\dd \samplealt | \sample)
    }]
    \\
    \leq &2\ex_{\sample \sim \prob}[
        \sum_{i \in I} \one_{\sample \in \mset^\Delta_i}
        \sum_{j \in J(i)}
        {
        {\coupling(\msetalt_j^\Delta | \sample)} \Delta
    }] = 2\Delta\,,
\end{align}
which vanishes as $\Delta$ goes to $0$. 
This proves point (2) of the result. 

Using the same decomposition of the integral, %
we get the same inequality for %
$\cost$,
\begin{align}
    \abs*{\int \cost \, \dd \coupling^\Delta - \int \cost \, \dd \coupling}
    \leq  \Delta\,,
\end{align}
which shows the convergence of the integral $\ex_{\coupling^\Delta}[\cost]$ to $\ex_{\coupling}[\cost]$, which is point (3) of the result.

Finally, take $\func \in \measfuncs$ \ac{lsc} and such that $\func_- \in \modspace[\cost]$, this means that there exists $\const > 0$ such that, $\base[\coupling]\; \als$,
$\func \geq - \const ( 1 + \base[\cost] + \cost)\,.$
But by construction, $\ex_{\coupling^\Delta}[1 + \base[\cost] + \cost] = 1 + \ex_{\prob}[\cost(\cdot, \base[\sample])] + \ex_{\coupling_\Regparam}[\cost]$ which goes to $\ex_{\coupling}[1 + \base[\cost] + \cost]$ as $\Delta \to 0$.
The last part of the result then follows from \citet[Lemma~4.3]{Villani2008OptimalTO} with $c \gets - \const(1 + \base[\cost] + \cost)$ and $h \gets \func$.
\end{proof}

\begin{lemma}\label{lemma:reformulate-dual}
Under the assumptions of \cref{theorem:app_weakened_general_regularized_duality}, the values of the two problems coincide:
	\begin{align}\label{eq:apxlem1}
		\sup \left\{
		\inner{\coupling}{\tildeobj} - \Reg(\coupling)
		:
		\coupling \in \wstarcl{\couplings_\prob \cap \measfuncs\cap\dspace}\,, \inner{\coupling}{\cost} + \Regalt(\coupling) \leq \radius
		\right\}
	\end{align}
    and	
    \begin{align}\label{eq:apxlem2}
		\sup \left\{
		\inner{\coupling}{\tildeobj} - \Reg(\coupling)
		:
		\coupling \in\couplings_\prob \,, \inner{\coupling}{\cost} + \Regalt(\coupling) \leq \radius
		\right\}.
\end{align}
\end{lemma}
\begin{proof} \leavevmode
We divide the proof in two parts, based on the two regularity assumption in~\cref{theorem:app_weakened_general_regularized_duality}. 

\noindent\underline{Case $(a)$:} \emph{$\Reg$ is $\weaktop{\dspace}{\contfuncs \cap \modspace[\cost]}$-\ac{usc}, $\obj$ is \ac{lsc} and $\tildeobj_- \in \modspace[\cost]$.} Thanks to the existence of a strictly feasible point $\coupling_\Regalt$ ($\inner{\coupling_\Regalt}{\cost} + \Regalt(\coupling_\Regalt) < \radius$) and the convexity (resp. concavity) of the constraint (resp. objective), \cref{lemma:ineq-constraints} ensures that the values of the two problems do not change when the inequality constraints are replaced by strict inequality constraints.

The value of \eqref{eq:apxlem1} is thus equal to
\begin{align}
        &\sup \left\{
		\inner{\coupling}{\tildeobj} - \Reg(\coupling)
		:
		\coupling \in \wstarcl{\couplings_\prob \cap \measfuncs\cap\dspace}\,, \inner{\coupling}{\cost} + \Regalt(\coupling) < \radius
		\right\}\\
        =
		&\sup \left\{
		\inner{\coupling}{\tildeobj} - \Reg(\coupling) - \indicator{\inner{\cdot}{\cost} + \Regalt < \radius}(\coupling)
		:
		\coupling \in \wstarcl{\couplings_\prob \cap \measfuncs\cap\dspace}
		\right\}\,,
\end{align}
with $\inner{\cdot}{\tildeobj} - \Reg - \indicator{\inner{\cdot}{\cost} + \Regalt < \radius}$ now $\weaktop{\dspace}{\pspace}$-\ac{lsc} since $\contfuncs \cap \modspace[\cost] \subset \pspace$ so that $\weaktop{\dspace}{\contfuncs \cap \modspace[\cost]}$ is weaker than $\weaktop{\dspace}{\pspace}$.
Thus, the weak-$*$ closure can be removed, and \eqref{eq:apxlem1} is equal to
\begin{align}
		&\sup \left\{
		\inner{\coupling}{\tildeobj} - \Reg(\coupling) - \indicator{\inner{\cdot}{\cost} + \Regalt < \radius}(\coupling)
		:
		\coupling \in \couplings_\prob \cap \measfuncs\cap\dspace
		\right\}\\
        =
		&\sup \left\{
		\inner{\coupling}{\tildeobj} - \Reg(\coupling)
		:
		\coupling \in\couplings_\prob \cap \measfuncs\cap\dspace \,, \inner{\coupling}{\cost} + \Regalt(\coupling) < \radius
		\right\}
		\\
        =
		&\sup \left\{
		\inner{\coupling}{\tildeobj} - \Reg(\coupling)
		:
		\coupling \in\couplings_\prob \cap \measfuncs \,, \inner{\coupling}{\cost} + \Regalt(\coupling) < \radius
		\right\},
\end{align}
where the last equality follows from \cref{lemma:criterion-dual-space} since $ \Regalt(\coupling)$ and  $\inner{\coupling}{\cost} $ have to be both finite for $\coupling$ to be feasible and $\Reg\propto\Regalt$.
All that is left to show is the equality
    \begin{align}
        &\sup \left\{
		\inner{\coupling}{\tildeobj} - \Reg(\coupling)
		:
		\coupling \in\couplings_\prob \cap \measfuncs \,, \inner{\coupling}{\cost} + \Regalt(\coupling) < \radius
		\right\} \\
        =
		&\sup \left\{
		\inner{\coupling}{\tildeobj} - \Reg(\coupling)
		:
		\coupling \in\couplings_\prob \,, \inner{\coupling}{\cost} + \Regalt(\coupling) < \radius
		\right\}\,.
\end{align}

The inequality $(\leq )$ is straightforward. Concerning the other side, for any $\coupling \in \couplings_\prob$ such that $\inner{\coupling}{\cost} + \Regalt(\coupling) < \radius$, we can consider $(\coupling_\run)_\run$ the sequence of probability distributions in $\couplings_\prob \cap \measfuncs$ given by \cref{lemma:dense}.
Thanks to $\cost$ belonging to $\contfuncs \cap \modspace[\cost]$ and the \ac{usc} assumption on $\Regalt\propto\Reg$, $\limsup_{\run \to \infty}  \inner{\coupling_\run}{\cost} + \Regalt(\coupling_\run) \leq 
 \inner{\coupling}{\cost} + \Regalt(\coupling) < \radius$ showing that $\coupling_\run$ is strictly feasible after some time.
 Moreover, since $\tildeobj$ is \ac{lsc} and $\tildeobj_-$ belongs to $\modspace[\cost]$, we get that
\begin{align}
		\inner{\coupling}{\tildeobj} - \Reg(\coupling)
        &\leq 
        \liminf_{\run \to \infty}~
		\inner{\coupling_\run}{\tildeobj} - \Reg(\coupling_\run)\\
        &\leq 
        \sup \left\{
		\inner{\coupling}{\tildeobj} - \Reg(\coupling)
		:
		\coupling \in\couplings_\prob \cap \measfuncs \,, \inner{\coupling}{\cost} + \Regalt(\coupling) < \radius
		\right\}\,,
\end{align}
which achieves the proof.

\medskip

\noindent\underline{Case $(2)$:} \emph{$\dom \Reg \subset \measfuncs$ and, for any sequence $\func_\run \in \dom \modular_\Reg$ non-negative such that $\func_\run \to 0$ as $\run \to \infty$, $\modular_\Reg(\func_\run) \to 0$.} In this case, $\modular_\Reg$ satisfies (H4) of \citet[Thm.~22]{pennanen2022topological} so that $(\heartmodspace[\Reg])^* = \kdmodspace[\Reg]$. As a consequence, $\couplings_\prob \cap \measfuncs\cap\dspace$ is $\weaktop{\dspace}{\pspace}$ closed so that $$ \wstarcl{\couplings_\prob \cap \measfuncs\cap\dspace} =  {\couplings_\prob \cap \measfuncs\cap\dspace}.$$
Therefore, using \cref{lemma:criterion-dual-space}, one gets that 
\begin{align}
		\sup \left\{
		\inner{\coupling}{\tildeobj} - \Reg(\coupling)
		:
		\coupling \in \couplings_\prob \cap \measfuncs\,, \inner{\coupling}{\cost} + \Regalt(\coupling) \leq \radius
		\right\}\,,
	\end{align}
and the result follows from the assumption on the domain of $\Reg$.
\end{proof}

\subsection{Examples}
\label{apx:examples}

We finally develop %
the two examples %
of \cref{sec:ex} by choosing carefully
the %
coupling 
$\base[\coupling]$.

\begin{corollary}[Cost-regularized WDRO]
   Let \cref{ass3} holds.
	Suppose that the function $\obj : \samples \to \R_+$ in $L^1(\prob)$, \ac{lsc}, and such that the function
            $(\sample, \samplealt) \mapsto {\obj(\samplealt)}/{\big(1 + \base[\cost](\sample, \samplealt) + \cost(\sample, \samplealt)\big)}$         
    is bounded, above and below, over $\samples^2$. 
    Then, for any $\reg, \regalt > 0$, the problem
	\begin{align}
		\sup \left\{
        \ex_{\coupling_2}[\obj]
		- \reg \ex_{\coupling}[\cost]
		:
		\coupling \in \couplings_\prob\,, \ex_{\coupling}[\cost] + \regalt\ex_{\coupling}[\cost] \leq \radius
		\right\}
	\end{align}
	coincides with its dual
	\begin{align}
		\inf \left\{
        \dualvar \radius +
		\ex_{\sample \sim \prob} \left[
            \esssup[\lebmeas]
        \left\{
			\obj - (\reg + (1+\regalt)\dualvar) \cost(\sample, \cdot)
        \right\}
        \right] :
		\dualvar \geq 0
		\right\}\,.
	\end{align}
\end{corollary}

\begin{proof}
In order to prove this result, we choose $\base \in \couplings_\prob$ such that, for any $\sample \in \samples$,
    \begin{align}
        \base(\dd \samplealt | \sample) \propto \frac{\one_{\samples}(\samplealt)}{1 + \base[\cost](\sample, \samplealt) + \cost(\sample, \samplealt)}e^{-\norm{\samplealt}^2/2}\dd \samplealt\,.
    \end{align}
We also define $\Reg$ for all $\coupling\in (\lebspb)^*$,
\begin{align}
    \Reg(\coupling) = \begin{cases}
        \reg\inner{\coupling}{\cost} &\text{if } \coupling \in \dmodspace[\cost] \text{ non-negative}\\
        \reg\int \cost \, \dd \coupling &\text{if } \coupling \in \measures \text{ non-negative}\\
        + \infty &\text{otherwise}\,.
    \end{cases}
\end{align}
With this definition, $\Reg$ satisfies the assumption of \cref{lemma:reg-modular-space} with $\coupling_\Reg = 0$ and we have
\begin{align}
    \modular_\Reg(\func) = \sup_{\coupling \in \measfuncs: \coupling \geq 0} \int (\abs{\func} - \reg \cost) \coupling \dd \base = \indicator{\abs \func \leq  \reg \cost}\,,
\end{align}
and thus $\heartmodspace[\Reg] = \{0\}$ so that $\pspace$ is equal to $\modspace[\cost]$.

It remains to check the regularity conditions on $\Reg$. $\Reg\restriction_{\dspace}$ is $\weaktop{\dspace}{\pspace}$-\ac{lsc} since $\cost$ belongs to $\pspace$ and the set of non-negative linear forms in $\dspace$ is $\weaktop{\dspace}{\pspace}$-closed.
On $\couplings_\prob$, $\Reg$ coincides with $\coupling \mapsto \reg \int \cost \, \dd \coupling$ and $\cost$ is inside $\contfuncs \cap \modspace[\cost]$ so it is actually $\weaktop{\dspace}{\contfuncs \cap \modspace[\cost]}$ continuous and thus regularity condition (a) of \cref{theorem:app_weakened_general_regularized_duality} holds. Finally, the strict feasibility condition is satisfied thanks to the distribution $\coupling(\dd \sample, \dd \samplealt) = \prob(\dd\sample)\dirac{\sample}(\dd \samplealt)$ so \cref{theorem:app_weakened_general_regularized_duality} applies.
Since $\lebsp[\infty] \subset \dspace$, the preconjugates of $\Reg + \dualvar \Regalt$ is given by,
\begin{align}
    (\Reg \restriction_{\dspace} + \dualvar \Regalt \restriction_{\dspace})_*(\func) = \indicator{\func \leq  (\reg+\dualvar \regalt) \cost}\,.
\end{align}
Hence, since the function $(\reg + \dualvar \regalt)\cost$ belongs to $\pspace = \modspace[\cost]$,
\begin{align}
    &\inf \left\{
        \dualvar \radius +
		\ex_{\sample \sim \prob} \left[
            \esssup
        \left\{
			\obj - \dualvar \cost(\sample, \cdot) - \func(\sample, \cdot)
        \right\}
        \right]
        +
        \left(\Reg\restriction_{\dspace} + \dualvar \Regalt\restriction_{\dspace}\right)_*(\func)
        :
		\func \in \pspace
		\right\}\\
&=
	\ex_{\sample \sim \prob} \left[
            \esssup[\lebmeas]
        \left\{
			\obj - (\reg + (1+\regalt)\dualvar) \cost(\sample, \cdot)
        \right\}
        \right]\,,
\end{align}
which gives the result.
\end{proof}

\begin{corollary}[$\phi$-divergence-regularized \ac{WDRO}]
    Let \cref{ass3} hold and take $\regparam,\regparamalt > 0$. Consider $\base \in \couplings$ with $\ex_{\base}[\cost]< +\infty$, a convex \ac{lsc} function $\phi : \R \to \extR$ with such that $\phi(1) = 0$, $[1, +\infty) \subset \dom \phi \subset [0, +\infty)$, $\phi'(+ \infty) = + \infty$ and define the associated divergence %
    \begin{align}
        \forall \coupling \in \couplings, \qquad \phidiv{\coupling}{\base} \defeq 
        \begin{cases}
            \int_{\samples^2} \phi \parens*{\frac{\dd \coupling }{\dd \base}} \dd \base &\text{if $\coupling$ is absolutely continuous \wrt $\base$}\\
            +\infty &\text{otherwise}\,.
        \end{cases}
    \end{align}
    Assume also that $\obj$ satisfies the growth condition:
    \begin{align}
        \exists \func \in \modspace[\cost]\,, \forall \alpha > 0\,, \int_{\samples^2} \phi^*\parens*{\alpha\abs{\tildeobj - \func}} \dd \base < + \infty\,.
    \end{align}
    Then, with $\Reg\gets\!\reg\phidiv{\,\cdot\,}{\base}$ and $\Regalt\gets\!\regalt\phidiv{\,\cdot\,}{\base}$, if\;$\eqref{eq:primal}$\;is\;strictly\;feasible,\;its\;value\;is\;equal\;to 
    \begin{align}
    \inf_{\dualvar \geq 0} \inf_{\bivarfuncalt \in \contfuncs} \dualvar \radius + \ex_{\sample \sim \prob} [\sup_{\altsample \in \samples} \obj(\altsample) - \dualvar \cost(\sample, \altsample) - \bivarfuncalt(\sample, \altsample)] + (\regparam + \dualvar \regparamalt)\int_{\samples^2} \phi^* \parens*{\frac{\bivarfuncalt(\sample, \samplealt)}{\regparam + \dualvar \regparamalt}} \dd \base (\sample, \samplealt)\,.
    \end{align}
\end{corollary}

\begin{proof}
    $D_\phi$ admits the following variational formula, see \eg \citet[Prop.~4.2.8]{agrawal2021optimal} applied with the pair $(\lebsp[\infty], \measures)$,
    \begin{align}
        \forall \coupling \in \measures,\quad 
        \phidiv{\coupling}{\base} = \sup \setdef*{\inner{\coupling}{\bivarfuncalt}  - \int_{\samples^2} \phi^* \circ \bivarfuncalt ~ \dd \base}{\bivarfuncalt \in \lebsp[\infty]}\,.
    \end{align}
    As seen in the second case of \cref{lemma:reformulate-dual}, $\dom\Reg \subset \measfuncs$ so that $\dspace \subset \measfuncs$. Hence, this variational formula applies to the whole $\dspace$ and since $\lebsp[\infty] \subset \modspace[\cost] \subset \pspace$, $\phidiv{\cdot}{\base}$ is $\weaktop{\dspace}{\pspace}$-\ac{lsc}.
    Moreover, \citet[Prop.~4.2.6]{agrawal2021optimal} with the %
    $(\pspace, \measures)$ and $(\pspace, \dspace)$ gives %
    \begin{align}
        \modular_\Reg(\func) = (\Reg \restriction_{\dspace})_*(\abs{\func}) =  \int_{\samples^2} \phi^* \circ \func ~ \dd \base\,.
    \end{align}
    The continuity %
    at 0 of $\modular_\Reg$ follows from the dominated convergence theorem and the fact that\;$\phi^*$ is non-decreasing.
    The growth condition %
    means that $\tildeobj$ lies in $\pspace$, and we can apply \cref{theorem:app_weakened_general_regularized_duality}.
\end{proof}
\end{document}